\documentclass[12pt]{article}
\usepackage[utf8]{inputenc}
\usepackage{color}
\definecolor{note_fontcolor}{rgb}{0.800781, 0.800781, 0.800781}
\usepackage{amsmath}
\usepackage{amsthm}
\usepackage{amssymb}
\usepackage[letterpaper]{geometry}
\PassOptionsToPackage{normalem}{ulem}
\usepackage{ulem}
\usepackage[unicode=true,pdfusetitle, bookmarks=true,bookmarksnumbered=false,bookmarksopen=false,
 breaklinks=false,pdfborder={0 0 1},backref=page,colorlinks=true]
 {hyperref}
\hypersetup{
 linkcolor=blue, citecolor=blue}

\usepackage[colorinlistoftodos]{todonotes} 
\usepackage{ stmaryrd }

\makeatletter

\numberwithin{equation}{section}
\numberwithin{figure}{section}

\theoremstyle{plain}
\newtheorem{thm}{\protect\theoremname}[section]
\theoremstyle{plain}
\newtheorem{lem}[thm]{\protect\lemmaname}
\theoremstyle{plain}
\newtheorem{cor}[thm]{\protect\corollaryname}
\theoremstyle{remark}
\newtheorem{rem}[thm]{\protect\remarkname}
\theoremstyle{plain}
\newtheorem{prop}[thm]{\protect\propositionname}
\theoremstyle{definition}
\newtheorem{defn}[thm]{\protect\definitionname}
\theoremstyle{plain}
\newtheorem{fact}[thm]{\protect\factname}
\newtheorem{assumption}[thm]{Assumption}

\@ifundefined{date}{}{\date{}}
\usepackage{graphicx}
\usepackage{appendix}

\usepackage{enumitem}
\usepackage[backref=page]{hyperref}
\usepackage{mathabx}

\usepackage{ytableau}

\usepackage{multirow}

\usepackage{appendix}

\usepackage{bbm}
\usepackage{xcolor}
\usepackage{tikz}
\usetikzlibrary{calc}
\usetikzlibrary{decorations.markings}
\usetikzlibrary{arrows.meta}
\usetikzlibrary{bending}

\def\centerarc[#1](#2)(#3:#4:#5)%
{ \draw[#1] ($(#2)+({#5*cos(#3)},{#5*sin(#3)})$) arc (#3:#4:#5); }

\tikzset{myptr/.style={decoration={markings,mark=at position 0.55 with %
{\arrow[scale=1,>=latex]{>}}},postaction={decorate}}}

\tikzset{myptr2/.style={decoration={markings,mark=at position 0.8 with %
{\arrow[rotate=-20,yshift=-.5\pgflinewidth,scale=0.75,>=latex]{>}}},postaction={decorate}}}

\tikzset{myptr2b/.style={decoration={markings,mark=at position 0.8 with %
{\arrow[rotate=20,yshift=.5\pgflinewidth,scale=0.75,>=latex]{>}}},postaction={decorate}}}

\newcommand{\hgline}[2]{
\pgfmathsetmacro{\thetaone}{mod(#1,360)}
\pgfmathsetmacro{\thetatwo}{mod(#2,360)}
\pgfmathsetmacro{\theta}{(\thetaone+\thetatwo)/2}
\pgfmathsetmacro{\phi}{abs(\thetaone-\thetatwo)/2}
\pgfmathsetmacro{\close}{less(abs(\phi-90),0.0001)}
\ifdim \close pt = 1pt
    \draw[blue!50] (\thetaone:1) -- (\thetatwo:1);
\else
    \pgfmathsetmacro{\R}{tan(\phi)}
    \ifdim \R pt < 0pt
        \pgfmathsetmacro{\distance}{sqrt(1+\R*\R)}
        \draw[blue!50] (\theta:-\distance) circle (\R);
    \else \ifdim \R pt > 0pt
        \pgfmathsetmacro{\distance}{sqrt(1+\R^2)}
        \draw[blue!50] (\theta:\distance) circle (\R);
        \fi
    \fi
\fi
}

\usepackage[font=small, margin=.5cm]{caption}

\usepackage{setspace}
\setdisplayskipstretch{1.5}
\makeatother

\providecommand{\corollaryname}{Corollary}
\providecommand{\definitionname}{Definition}
\providecommand{\factname}{Fact}
\providecommand{\lemmaname}{Lemma}
\providecommand{\propositionname}{Proposition}
\providecommand{\remarkname}{Remark}
\providecommand{\theoremname}{Theorem}

\begin{document}
\global\long\def\F{\mathrm{\mathbf{F}} }%
\global\long\def\Aut{\mathrm{Aut}}%
\global\long\def\C{\mathbf{C}}%
\global\long\def\H{\mathcal{H}}%
 
\global\long\def\U{\mathcal{U}}%

\global\long\def\ext{\mathrm{ext}}%
 
\global\long\def\hull{\mathrm{hull}}%
 
\global\long\def\triv{\mathrm{triv}}%
 
\global\long\def\Hom{\mathrm{Hom}}%

\global\long\def\trace{\mathop{\mathrm{tr}}}%
\global\long\def\End{\mathrm{End}}%
 
\global\long\def\tsg{\widetilde{\Sigma_{g}}}%

\global\long\def\L{\mathcal{L}}%
\global\long\def\W{\mathcal{W}}%
\global\long\def\E{\mathbb{E}}%
\global\long\def\SL{\mathrm{SL}}%
\global\long\def\R{\mathbf{R}}%
\global\long\def\Pairs{\mathrm{PowerPairs}}%
\global\long\def\Z{\mathbf{Z}}%
\global\long\def\rs{\to}%
\global\long\def\A{\mathcal{A}}%
\global\long\def\a{\mathsf{a}}%
\global\long\def\rsa{\rightsquigarrow}%

\global\long\def\D{\mathcal{D}}%
\global\long\def\b{\mathbf{b}}%
\global\long\def\df{\mathrm{def}}%
\global\long\def\eqdf{\stackrel{\df}{=}}%
\global\long\def\ZZ{\overline{Z}}%
\global\long\def\Tr{\mathop{\mathrm{Tr}}}%
\global\long\def\N{\mathbf{N}}%
\global\long\def\std{\mathrm{std}}%
 
\global\long\def\HS{\mathrm{H.S.}}%
\global\long\def\e{\mathbf{e}}%
\global\long\def\c{\mathbf{c}}%
\global\long\def\d{\mathbf{d}}%
\global\long\def\AA{\mathbf{A}}%
\global\long\def\BB{\mathbf{B}}%
 
\global\long\def\v{\mathbf{v}}%
\global\long\def\spec{\mathrm{spec}}%
\global\long\def\Ind{\mathrm{Ind}}%
\global\long\def\half{\frac{1}{2}}%
\global\long\def\Re{\mathrm{Re}}%
 
\global\long\def\Im{\mathrm{Im}}%
\global\long\def\Rect{\mathrm{Rect}}%
\global\long\def\Crit{\mathrm{Crit}}%
\global\long\def\Stab{\mathrm{Stab}}%
\global\long\def\SL{\mathrm{SL}}%
\global\long\def\Tab{\mathrm{Tab}}%
\global\long\def\Cont{\mathrm{Cont}}%
\global\long\def\I{\mathcal{I}}%
\global\long\def\J{\mathcal{J}}%
\global\long\def\short{\mathrm{short}}%
\global\long\def\Id{\mathrm{Id}}%
\global\long\def\B{\mathcal{B}}%
\global\long\def\ax{\mathrm{ax}}%
\global\long\def\cox{\mathrm{cox}}%
\global\long\def\row{\mathrm{row}}%
\global\long\def\col{\mathrm{col}}%
\global\long\def\X{\mathbb{X}}%
\global\long\def\Fat{\mathsf{Fat}}%

\global\long\def\V{\mathcal{V}}%
\global\long\def\P{\mathbb{P}}%
\global\long\def\Fill{\mathsf{Fill}}%
\global\long\def\fix{\mathsf{fix}}%
 
\global\long\def\reg{\mathrm{reg}}%
\global\long\def\edge{E}%
\global\long\def\id{\mathrm{id}}%
\global\long\def\emb{\mathrm{emb}}%

\global\long\def\Hom{\mathrm{Hom}}%
 
\global\long\def\F{\mathrm{\mathbf{F}} }%
  
\global\long\def\pr{\mathrm{Prob} }%
 
\global\long\def\tr{\mathop{\mathrm{Tr}}}%
\global\long\def\core{\mathrm{Core}}%
\global\long\def\pcore{\mathrm{PCore}}%
\global\long\def\im{\vartheta}%
\global\long\def\br{\mathsf{BR}}%
 
\global\long\def\sbr{\mathsf{SBR}}%
 
\global\long\def\ebs{\mathsf{EBS}}%
 
\global\long\def\ev{\mathrm{ev}}%
 
\global\long\def\CC{\mathcal{C}}%
 
\global\long\def\sides{\mathrm{Sides}}%
\global\long\def\tp{\mathrm{top}}%
\global\long\def\lf{\mathrm{left}}%
\global\long\def\MCG{\mathrm{MCG}}%
\global\long\def\EE{\mathcal{E}}%
 
\global\long\def\mog{\mathfrak{MOG}}%
 
\global\long\def\fg{\le_{\mathrm{f.g.}}}%
 
\global\long\def\v{\mathfrak{v}}%
\global\long\def\e{\mathfrak{e}}%
\global\long\def\f{\mathfrak{f}}%
 
\global\long\def\d{\mathfrak{d}}%
\global\long\def\he{\mathfrak{he}}%

\global\long\def\defect{\mathrm{Defect}}%
\global\long\def\M{\mathcal{M}}%
\global\long\def\sdefect{\max\defect}%
\global\long\def\p{\mathfrak{p}}%
\global\long\def\free{\mathbf{F}}%
\global\long\def\Q{\mathbf{Q}}%
\global\long\def\Part{\mathrm{Part}}%
\global\long\def\subperm{\mathrm{SubPerm}}%
\global\long\def\ntr{\mathop{\mathrm{tr}}}%

\global\long\def\F{\mathbb{\mathbb{\mathbf{F}}}}%
 
\global\long\def\defi{\stackrel{\text{def}}{=}}%
 
\global\long\def\G{\Gamma}%
 
\global\long\def\cay{\text{Cay}(\G)}%
\global\long\def\H{\mathbb{H}^{2}}%
\global\long\def\supp{\mathop{\text{supp}}}%

\title{Strong convergence of uniformly random\\
permutation representations of surface groups}

\author{Michael Magee, Doron Puder, and Ramon van Handel}

\maketitle

\begin{abstract}
Let $\Gamma$ be the fundamental group of a closed orientable surface
of genus at least two. Consider the composition of a uniformly random element of $\Hom(\Gamma,S_n)$ with the $(n-1)$-dimensional irreducible representation of $S_n$. We prove the strong convergence in probability as $n\to\infty$ of this sequence of random representations to the regular representation of $\Gamma$.

As a consequence, for any closed orientable hyperbolic surface $X$, with probability tending to one as $n\to\infty$, a uniformly random degree-$n$ covering space of $X$ has near optimal relative spectral gap --- ignoring the eigenvalues that arise from the base surface $X$.

To do so, we show that the polynomial method of proving strong convergence can be extended beyond rational settings.

To meet the requirements of this extension we prove two new kinds
of results. First, we show there are effective polynomial approximations of expected values of traces of elements of $\Gamma$ under random homomorphisms to $S_{n}$. Secondly, we estimate the growth rates of probabilities that a finitely supported random walk on $\Gamma$ is a proper power after a given number of steps.

\tableofcontents{}
\end{abstract}

%
%

\section{Introduction}

In the following $\Gamma=\Gamma_{g}$ will be the fundamental group of a closed orientable surface of genus $g\geq2$. We view $g$ and $\Gamma$ as fixed in the paper. We denote by $S_{n}$ the group
of permutations of $[n]\eqdf\{1,\ldots,n\}$ and $\U(n)$ the group of $n\times n$ complex unitary matrices. The $n$-dimensional linear representation of $S_{n}$ by permutation matrices is not irreducible but has an $(n-1)$-dimensional irreducible subrepresentation 
\[
\std:S_{n}\to\U(n-1)
\]
obtained by removing the non-zero invariant vectors. The group $\Gamma$ also has an important unitary representation called the \emph{regular
representation}:
\[
\lambda:\Gamma\to\U(\ell^{2}(\Gamma)),\quad\lambda(\gamma)[f](x)\eqdf f(\gamma^{-1}x).
\]
For $n\in\N$ the set of homomorphisms 
\[
\X_{g,n}\eqdf\Hom(\Gamma_{},S_{n})
\]
is finite and we equip it with the uniform probability measure. Let $\C[\Gamma]$ denote the group algebra of $\Gamma$. Every linear
representation of $\Gamma$ extends linearly to one of $\C[\Gamma]$.

\begin{thm} \label{thm:strong_convergence}
For all $x\in\C[\Gamma]$, for $\phi_{n}\in\Hom(\Gamma,S_{n})$
uniformly random, as $n\to\infty$
\[
\|\std\circ\phi_{n}(x)\|\to\|\lambda(x)\|
\]
in probability. The norms on each side are operator norms.
\end{thm}

In the language of strong convergence of unitary representations,
Theorem \ref{thm:strong_convergence} says precisely that the random
representations $\std\circ\phi_{n}$ converge strongly in probability
to $\lambda$. See \cite{magee2025strongconvergenceunitarypermutation}
for a survey on strong convergence. Theorem \ref{thm:strong_convergence}
is the analog, for surface groups, of the main theorem of Bordenave--Collins
\cite{BordenaveCollins} for free groups. In \S\ref{subsec:Prior-results}
we cover the history of strong convergence results prior
to Theorem \ref{thm:strong_convergence}. 

While we state Theorem \ref{thm:strong_convergence} here in a qualitative form for simplicity, the proof yields a quantitative result with a polynomial convergence rate, see Theorem \ref{thm:strongquant}.

\medskip

\emph{Random hyperbolic surfaces.} One of the primary motivations for Theorem \ref{thm:strong_convergence} is a known important consequence on the spectral theory of random hyperbolic surfaces. A hyperbolic surface $X$ is a complete Riemannian surface of constant curvature $-1$. Such a surface has an associated Laplacian operator $\Delta_{X}$.
If $X$ is closed, the spectrum of $\Delta_{X}$ as an unbounded operator on $L^{2}(X)$ is discrete consisting of real eigenvalues that accumulate only at infinity. If $X$ is connected, the bottom of the spectrum of $\Delta_X$ is a simple eigenvalue at $0$ of multiplicity one.

The first non-zero eigenvalue, $\lambda_1(X)$, is of fundamental significance to e.g.~the dynamics of the geodesic flow on $X$, the statistics of lengths of closed geodesics, and notions of connectedness of $X$. In cases where $X$ is a congruence hyperbolic surface, Selberg recognized \cite{zbMATH03229797} that the condition $\lambda_1(X)\geq\frac{1}{4}$
--- the Selberg Eigenvalue conjecture --- is the archimedean analog of the Ramanujan--Petersson conjectures from number theory, and from the point of view of the Selberg zeta function is also the analog of the Riemann hypothesis for these surfaces.

If $\{X_n\}_{n=1}^{\infty}$ are a sequence of closed hyperbolic
surfaces with genera tending to $\infty$, then $\limsup_{n\to\infty}\lambda_{1}(X_{n})\leq\frac{1}{4}$ \cite{Huber} and hence $\frac{1}{4}$ is the asymptotically optimal value of $\lambda_{1}(X_{n})$ for such surfaces. The number $\frac{1}{4}$ arises here as the bottom of the spectrum of the Laplacian on hyperbolic space of dimension two.

\begin{thm} \label{thm:main-hyperbolic}
Let $X$ be a closed connected orientable hyperbolic surface. Let $X_n$ denote a uniformly random degree $n$ covering space of $X$.\footnote{If $\pi_1(X)\cong \Gamma$, then by definition, a uniformly random degree-$n$ cover of $X$ is given by a uniformly random element of $\X_{g,n}$ via the bijection between $\X_{g,n}$ and covers with a labeled fiber; see \cite[\S1]{MPasympcover}.} With probability tending to one as $n\to\infty$
\[
\spec(\Delta_{X_{n}})\cap\left[0,\frac{1}{4}-o(1)\right)=\spec(\Delta_{X})\cap\left[0,\frac{1}{4}-o(1)\right),
\]
with the same multiplicities on both sides. \\
In particular, if $\lambda_{1}(X)\geq\frac{1}{4}$, then $\lambda_{1}(X_{n})\to\frac{1}{4}$
in probability as $n\to\infty$.
\end{thm}

We omit the proof that Theorem \ref{thm:main-hyperbolic} follows from Theorem \ref{thm:strong_convergence}, which follows as in \cite[Appendix A]{zbMATH07974824}.

The existence of covering spaces of hyperbolic surfaces that have optimal spectral gaps was proved in \cite{zbMATH07734891,zbMATH07974824}. What Theorem \ref{thm:main-hyperbolic} shows is that not only do such covers exist, but that this is in fact the \emph{typical} behavior of covering spaces of closed hyperbolic surfaces.
In \S\ref{subsec:Prior-results} we describe prior results to Theorem \ref{thm:main-hyperbolic}.

\medskip

\emph{Further developments.}
We highlight two further important developments that arise from the results of this paper. By a recent result of Hide--Moy--Naud \cite[Thm.~1.3]{hidemoynaud2025}, 
Theorem~\ref{thm:strong_convergence} also yields the analog of Theorem~\ref{thm:main-hyperbolic} for
surfaces $X$ with variable negative curvature. On the other hand, Hide--Macera--Thomas \cite{hidethomasmacera2025} apply Theorem~\ref{thm:strongquant} below to obtain a quantitative
form of Theorem \ref{thm:main-hyperbolic} that achieves the first polynomial convergence rate for
the optimal spectral gap of random hyperbolic surfaces. The significance of these
developments is discussed further in \S\ref{subsec:Prior-results}.

\subsection{Technical innovation I: The generalized polynomial method\label{subsec:Technical-innovation-I:surpassing-poly}}

The polynomial method, which was introduced in the recent work of Chen, Garza-Vargas, Tropp, and the third named author \cite{chen2024new,chen2024newapproachstrongconvergence}
and further developed by de la Salle and the first named author \cite{MageeSalle2}, provides a powerful approach for establishing strong convergence in situations that are outside the reach of previous methods. To date, however, all applications of the polynomial method have relied fundamentally on the property that the spectral statistics of the random matrices that arise in these applications are rational functions of $\frac{1}{n}$. This restrictive property fails manifestly for the model considered in this paper.\footnote{Indeed, in the setting of this paper, the spectral statistics are not even analytic as a function of $\frac{1}{n}$.
This can be verified by checking, by an explicit computation, that the coefficients $\a_i(\gamma)$ in Theorem~\ref{thm:MP-Laurent} grow faster than exponentially in $i$ when $\gamma$ is one of the standard generators of $\Gamma$. We omit the details.}

The basis for the present paper is a refinement of the polynomial method, which makes it applicable to much more general situations where the dimension dependence of the spectral statistics need not even be analytic. This greatly expands the potential range of applications of this method. More precisely, this method will give rise to the following general criterion for strong convergence.

Let $\Lambda$ be a finitely generated group with a finite generating set $S$. For $\gamma\in\Lambda$, denote by $|\gamma|$ the word length of $\gamma$ with respect to $S$. Denote by $\lambda$ the left regular representation of $\Lambda$, and for each $n\in\N$ let 
\[
    \pi_{n}:\Lambda\to \U(n)
\]
be a random unitary representation. We denote by $\tr M$ the trace of a matrix $M$ and by $\ntr M$ its normalized trace. We let $\tau$ be the canonical trace on the reduced group $C^*$-algebra $C_{{\rm red}}^{*}(\Lambda)$. (The definitions of these notions are recalled in \S\ref{sec:notation}.)

To establish that $\pi_{n}\to\lambda$ strongly in probability, we will need two assumptions. 
The first is an effective $\frac{1}{n}$-expansion for the expected spectral statistics.

\begin{assumption}
\label{ass:cond1}
There exists an absolute constant $C\ge1$, and
constants $u_k(\gamma)\in\R$ for all $k\in\Z_+$ and $\gamma\in\Lambda$, such that
$u_0(\gamma)=\tau(\lambda(\gamma))$ and
\[
        \Bigg|
        \E\big[{\ntr \pi_n(\gamma)}\big] -
        \sum_{k=0}^{q-1} \frac{u_k(\gamma)}{n^k}\Bigg|
        \le \frac{(Cq)^{Cq}}{n^q}
\]
for all $\gamma\in\Lambda$, $q\ge |\gamma|$, and $n\ge Cq^C$.
\end{assumption}

In the following, we will extend $u_k$ linearly to the group algebra
$\C[\Lambda]$. 
The second assumption is that $u_1$ is ``tempered'' in the sense of \cite{chen2024new,MageeSalle2}.

\begin{assumption}
\label{ass:cond2}
For every self-adjoint $x\in\C[\Lambda]$, we have
\[
        \limsup_{p\to\infty} |u_1(x^p)|^{1/p} \le
        \|\lambda(x)\|.
\]
\end{assumption}

The following theorem is proved in \S\ref{sec:Proof-of-new-criterion-strong-convergence}.

\begin{thm}
\label{thm:main} Suppose that Assumptions \ref{ass:cond1} and \ref{ass:cond2}
hold. Then 
\[
\|\pi_{n}(x)\|\le\|\lambda(x)\|+o(1)\quad\text{with probability}\quad1-o(1)
\]
as $n\to\infty$ for every $x\in\C[\Lambda]$. 
\end{thm}

Theorem \ref{thm:main} gives a strong convergence upper bound. If $C_{{\rm red}}^{*}(\Lambda)$ has a unique trace, this automatically implies the lower bound \cite[\S  5.3]{MageeSalle2}. 
We state the result here in a qualitative form for simplicity, but Theorem \ref{thm:quantmain} below in fact provides a quantitative form of this result with a polynomial convergence rate.

Let us highlight two novel features of the above result.
\begin{itemize}
\item
Previous applications of the polynomial method required that $\E[\ntr\pi_{n}(\gamma)]$ is a rational function of $\frac{1}{n}$. Assumption \ref{ass:cond1} is a drastic relaxation of this assumption: it requires only that the spectral statistics are in a Gevrey class (so they need not even be analytic, let alone rational) as a function of $\frac{1}{n}$. 
\item
To verify Assumption \ref{ass:cond1}, one only needs to consider individual group elements $\gamma\in\Lambda$ rather than elements of the group algebra $x=\sum_{\gamma}\alpha_{\gamma}\gamma\in\C[\Lambda]$.
\end{itemize}
The latter is crucial for applications. In principle, one may expect that a method that applies to polynomial (or rational) spectral statistics could be adapted to spectral statistics that are well approximated by polynomials. This is what Assumption~\ref{ass:cond1} ensures. However, a key feature of the polynomial method is that it captures cancellations between the terms of $x\in\C[\Lambda]$ in the $\frac{1}{n}$-expansion of $\E[\ntr\pi_{n}(x)]$ \cite[\S 2.3]{chen2024new}. It is not clear
that such cancellations are preserved when the spectral statistics are approximated by polynomials.
Since we do not know how to establish \emph{a priori} bounds that capture cancellations, it is
essential that this is not required by Assumption \ref{ass:cond1}.
We elaborate further on this point in the remarks at the end of \S\ref{sec:Proof-of-new-criterion-strong-convergence}.

\subsection{Technical innovation II: Effective polynomial approximation\label{subsec:Technical-innovation-II:Effective-approx}}

We aim to apply the criterion for strong convergence from Theorem \ref{thm:main} to the setting that $\Lambda=\Gamma=\Gamma_{g}$ is a surface group with $g\ge 2$, and 
\[
\pi_{n}=\mathrm{std}\circ\phi_{n},
\]
where $\phi_{n}$ is a uniformly distributed random element of $\mathrm{Hom}(\Gamma_{g},S_{n})$. Thus we need to establish Assumptions \ref{ass:cond1} and \ref{ass:cond2} for this setting.\footnote{There is a slight discrepancy between the present setting and that of \S\ref{subsec:Technical-innovation-I:surpassing-poly}, since here $\pi_n$ is of dimension $n-1$ rather than $n$. This issue is cosmetic in nature and will be dispensed with in \S\ref{sec:Main-proof}.} We begin with Assumption~\ref{ass:cond1}.

For $\sigma\in S_{n}$ let $\fix(\sigma)$ denote the number of fixed points of $\sigma$. For a random variable $f$ on $\X_{g,n}$ we denote by $\E_{g,n}[f]$ the expectation of $f$ with respect to the uniform measure. Natural random variables on $\X_{g,n}$ come from elements of $\gamma$: given $\gamma\in\Gamma$ we define
\begin{align*}
\fix_{\gamma} & :\X_{g,n}\to\Z\\
\fix_{\gamma} & :\phi\mapsto\fix(\phi(\gamma)).
\end{align*}
Then we have
\begin{equation}
\label{eq:fixrep}
\E[\tr \pi_{n}(\gamma)]=\E_{g,n}[\fix_{\gamma}]-1.
\end{equation}
An asymptotic expansion for $\E_{g,n}[\fix_{\gamma}]$
was obtained in the previous work \cite{MageePuderCore,MPasympcover}
of the first two authors.
\begin{thm}[{\cite[Thm.~1.1]{MPasympcover}}]
\label{thm:MP-Laurent}For each $\gamma\in\Gamma$ there exists
a (unique) sequence $\left\{ \a_{i}(\gamma)\right\} _{i=-1}^{\infty}$
such that for any $q\in\N$, as $n\to\infty$
\[
\E_{g,n}[\fix_{\gamma}]=\mathsf{a}_{-1}(\gamma)n+\a_{0}(\gamma)+\a_{1}(\gamma)n^{-1}+\cdots+\a_{q-1}(\gamma)n^{-(q-1)}+O_{q,\gamma}\left(n^{-q}\right),
\]
where the implied constant in the big O depends on $q$ and $\gamma$.
\end{thm}

In \S\ref{sec:EffectivePolynomialApproximation}, we prove the following theorem that strengthens Theorem \ref{thm:MP-Laurent}. Here and in the sequel, we fix a standard set of generators of $\Gamma$ of size $2g$ coming from the presentation
\begin{equation}
\Gamma=\left\langle a_{1},b_{1},\ldots,a_{g},b_{g}\,\mid\,[a_{1},b_{1}]\cdots[a_{g},b_{g}]\right\rangle ,\label{eq:presentation of Gamma}
\end{equation}
and for $\gamma\in\Gamma$ write $|\gamma|$ for the word length of $\gamma$ with respect to these generators. 

\begin{thm}
\label{thm:asymp-effective}There exists $C=C(g)>1$ such that for
all $q,n\in\N$ with $n\geq Cq^{C}$, and for all $\gamma\in\Gamma$ with
$|\gamma|\leq q$,
\[
\left|\E_{g,n}[\fix_{\gamma}]-\a_{-1}(\gamma)n-\a_{0}(\gamma)-\a_{1}(\gamma)n^{-1}-\cdots-\a_{q-1}(\gamma)n^{-(q-1)}\right|\leq(Cq)^{Cq}n^{-q}.
\]
\end{thm}

The proof of this effective asymptotic expansion proceeds in two parts.
The first part of the analysis relies on an exact expansion of $\E_{g,n}[\fix_{\gamma}]$
in terms of all irreducible representations of $S_{n}$ that is developed
in \cite{MPasympcover}. A careful analysis shows that we may truncate
this expansion only to those representations defined by Young diagrams
$\lambda\vdash n$ with $O(q)$ boxes outside the first row and column,
with an effective error bound as in Assumption \ref{ass:cond1}. This is
Proposition \ref{prop:tail-cut}.

It is not clear, however, how to analyze the truncated expansion using
the methods of \cite{MPasympcover}. Instead, we introduce a new method
for analyzing the truncated expansion by means of a combinatorial
integration formula for the symmetric group. This enables us to show
that the truncated expansion is in fact a rational function of $\frac{1}{n}$.
An effective error bound for the truncated expansion can be then achieved
using the analytic theory of polynomials, in a manner closely analogous
to the polynomial method for random matrix models with rational spectral
statistics \cite{chen2024new,chen2024newapproachstrongconvergence,MageeSalle2}.

\subsection{Technical innovation III: Geometry of proper powers\label{subsec:Technical-innovation-III:powers}}

It remains to establish Assumption \ref{ass:cond2} for surface groups.

For every element $1\ne\gamma\in\Gamma$, there is a unique $\gamma_0\in\Gamma$ and $k\in\N$ with $\gamma=\gamma_{0}^{\,k}$ such that $\gamma_{0}$ is not a proper power.\footnote{This follows immediately from the fact the non-trivial elements of $\Gamma$ admit cyclic centralizer. The latter fact can be found, e.g., in \cite[p.~23]{farb2011primer}.}
We denote by $\Gamma_{\rm np}$ the set of non-powers in $\Gamma$.
Let $\omega(k)$ be the number of positive divisors of $k$,
and by abuse of notation, we let
\[
    \omega(\gamma) \eqdf \begin{cases}
    0 & \text{if }\gamma=1,\\
    \omega(k) & \text{if }\gamma=\gamma_0^k\text{ with }\gamma_0\in\Gamma_{\rm np}.
    \end{cases}
\]
In the notation of Theorem \ref{thm:asymp-effective},
it is shown in \cite[Thm.~1.2]{MPasympcover} that $\a_{0}(\gamma)=\omega(\gamma)$.
Thus $u_1$ in Assumption \ref{ass:cond2} can be expressed 
for any $x=\sum_{\gamma\in\G}\alpha_{\gamma}\gamma\in\C[\Gamma]$ as
\begin{equation}
\label{eq:u1basic}
u_{1}(x)=\sum_{\gamma}\alpha_{\gamma}\left[\omega(\gamma)-1\right]
\end{equation}
using \eqref{eq:fixrep}, see \S\ref{sec:Main-proof}. Note that $u_1$ is determined only by the contributions of $\gamma\in\Gamma$ that are either the identity or a proper power, since $\omega(\gamma)-1=0$ for $\gamma\in\Gamma_{\rm np}$.

We prove the following theorem in \S\ref{sec:Proper-powers}. 

\begin{thm}
\label{thm:powers}For any self adjoint $x\in\C[\G]$, 
\[
\limsup_{p\to\infty}\left|u_{1}(x^{p})\right|^{1/p}\le\left\Vert \lambda(x)\right\Vert .
\]
\end{thm}

To prove Theorem \ref{thm:powers} it suffices, by \cite[Prop.~6.2]{MageeSalle2} and as surface groups satisfy the rapid decay property \cite[Thm.~3.2.1]{jolissaint1990rapidly},
to control the exponential
growth rate of the probability that a random walk in $\Gamma$ with an arbitrary finite set of generators lands on a proper power.

In the free group setting, the latter is readily accomplished by a
simple spectral argument that dates back to \cite[Lem.\ 2.4]{friedman2003relative}.
This argument relies on the following elementary property of free
groups: if a word $w=a_{i_{1}}\cdots a_{i_{q}}$ in the free generators
$a_{1},\ldots,a_{r}$ reduces to a proper power $\gamma_{0}^{k}$
and we write $\gamma_{0}=bhb^{-1}$ where $h$ is the cyclic reduction of $\gamma_0$,
then $w$ must be a concatenation of five consecutive sub-words that
reduce to $b,h,h,h^{k-2}$, and $b^{-1}$, respectively.

The difficulty in the setting of surface groups is that such a property is no longer true. However, we will show in Lemma \ref{lem:padding power words to get nice structure} that an approximate form of this property remains valid: if a word $w$ in the standard generators of $\Gamma$ is equivalent in $\G$ to a proper power $\gamma_{0}^{k}$, then there must exist $b,h\in\Gamma$ so that $w$ is a concatenation of five consecutive sub-words which are equivalent in $\G$ to elements that are $c\log|w|$-close to $b,h,h,h^{k-2}$, and $b^{-1}$ (here $u$ is $\ell$-close to $v$ if $u$ is equivalent in $\G$ to $\gamma v\gamma'$ with $|\gamma|,|\gamma'|\le\ell$). The proof of this fact makes crucial use of the hyperbolic geometry of $\mathrm{Cay}(\Gamma)$. This suffices to conclude the proof of Assumption \ref{ass:cond2} for surface groups.\footnote{After the first version of this paper appeared, it was shown by de la Salle \cite{delaSalle2025} that this approach can be generalized to the case that $\Gamma$ is any hyperbolic group.}

\subsection{Prior results\label{subsec:Prior-results}}

\subsubsection*{Strong convergence}

The first result in the spirit of Theorem \ref{thm:strong_convergence} was the
breakthrough result of Haagerup and Thorbj{\o}rnsen \cite{HaagerupThr},
who proved that i.i.d.\ GUE random matrices strongly converge
to free semicircular random variables. Let
$\F_{r}$ denote a free group of rank $r$. Collins and Male proved
in \cite{Collins2014} that Haar-distributed elements of $\Hom(\free_{r},\U(n))\cong\U(n)^{r}$
strongly converge to the regular representation of $\F_{r}$.
In another breakthrough \cite{BordenaveCollins}, Bordenave and Collins
proved that uniformly random elements of $\Hom(\F_{r},S_{n})$ composed
with $\std$ strongly converge to the regular representation
of $\F_{r}$.

Subsequent developments include results on
higher dimensional representations of $\U(n)$ or $S_{n}$ \cite{zbMATH07870147,chen2024new,MageeSalle2,cassidy2025randompermutationsactingktuples}, results related to Hayes' work \cite{zbMATH07565550} on the Peterson-Thom conjecture \cite{belinschi2024strongconvergencetensorproducts,bordenave2024normmatrixvaluedpolynomialsrandom,ParraudHayes,chen2024newapproachstrongconvergence}, and work of Austin on annealed almost periodic entropy \cite{austin2024notionsentropyunitaryrepresentations}.

Regarding the existence of finite dimensional unitary representations
that strongly converge to the regular representation, results
have been obtained for various classes of groups \cite{zbMATH07974824,zbMATH07929052,magee2023stronglyconvergentunitaryrepresentations}
other than free groups. Most relevant to the current paper is the
result of the first named author and Louder \cite[Cor.~1.2]{zbMATH07974824} that for surface groups $\Gamma$ as in this paper, there exists a sequence of $\{\phi_{n}\in\Hom(\Gamma,S_{n})\}_{n=1}^{\infty}$ such that $\std\circ\phi_{n}$ strongly converge to $\lambda$. Theorem \ref{thm:strong_convergence} shows not only that such homomorphisms exist, but that this is in fact the \emph{typical} behavior.

The recent works \cite{chen2024new,chen2024newapproachstrongconvergence}
introduced a new approach to strong convergence, called the \emph{polynomial
method}, that is based almost entirely on soft arguments, in contrast
to previous works that were based on problem-specific analytic methods
or delicate combinatorial estimates. This method has made it possible
to establish strong convergence in previously inaccessible situations
and to achieve improved quantitative results. These works form the basis
for the generalized polynomial method that is developed in the present paper (see \S\ref{subsec:Technical-innovation-I:surpassing-poly} and \S\ref{sec:Proof-of-new-criterion-strong-convergence}).

\subsubsection*{Hyperbolic surfaces}

Whether there exists a sequence of closed hyperbolic surfaces  $\{X_n\}_{n=1}^{\infty}$ with genera tending to $\infty$ such that $\lambda_{1}(X_{n})\to\frac{1}{4}$ was an old question of Buser \cite{Buser}. In view of analogous results on optimal spectral gaps of random regular graphs \cite{friedman2008proof,zbMATH07346527}, it is natural to approach this question through the study of random hyperbolic surfaces.
We presently discuss previous results on three models of hyperbolic surfaces:
\begin{enumerate}
\item
The Brooks--Makover (BM) model
and variations, obtained by constructing a random cusped hyperbolic surface out of gluing $m$ ideal triangles and then compactifying the cusps. 
\item The Weil--Petersson (WP) random model that comes from the Weil--Petersson Kähler form on the moduli space of genus $g$ hyperbolic surfaces.
\item Models of random degree $n$ covering spaces.
\end{enumerate}
We say events hold w.h.p.~(with high probability) if they hold with probability tending to one as the relevant parameter $(m/g/n)$ tends to infinity. In random cover models we understand that $\lambda_{1}$ is measured ignoring the eigenvalues of the base surface.

The first result that random closed hyperbolic surfaces in the BM
model have a uniform positive lower bound on $\lambda_{1}$
was due to Brooks--Makover \cite{zbMATH05033782}. The same result
for WP random surfaces was obtained by Mirzakhani in \cite{zbMATH06186878}.
The explicit bound $\lambda_{1}\geq\frac{3}{16}-o(1)$ w.h.p.\ for
closed surfaces was first obtained in \cite{MageeNaudPuder} for
the uniform random covering model. Following this, the same result was obtained
in the WP model by Wu--Xue \cite{zbMATH07511031} and Lipnowski--Wright
\cite{zbMATH07815241}.

The optimal bound $\lambda_{1}\geq\frac{1}{4}-o(1)$ w.h.p.\ was finally achieved by
Hide and the first named author \cite{zbMATH07734891} in a random cover
model for \emph{non-closed} hyperbolic surfaces. This suffices, by means of a compactification procedure, to settle Buser's conjecture in the affirmative.
The same proof applies verbatim to a model of gluing hyperbolic triangles\footnote{Although it was not pointed out at the time, the proof of \cite{zbMATH07734891} works just as well for random covers of the quotient of $\H$ by reflections in the sides of an ideal hyperbolic triangle. In this case, the input is the result of Bordenave and Collins for random homomorphisms from the free product $\Z/2\Z * \Z/2\Z * \Z/2\Z$ to $S_n$. This is not exactly the same as the Brooks--Makover gluing, however.} and as a result, one obtains $\lambda_{1}\geq\frac{1}{4}-o(1)$ w.h.p.~in close cousins of the BM model for closed surfaces. 
In \cite{zbMATH07974824},
the existence of covers of closed hyperbolic surfaces such that $\lambda_{1}\geq\frac{1}{4}-o(1)$
was proved without a compactification procedure. 

Very recently, it was shown in the impressive works of Anantharaman--Monk 
\cite{anantharaman2023friedmanramanujanfunctionsrandomhyperbolic,anantharaman2025friedmanramanujanfunctionsrandomhyperbolic}
that $\lambda_{1}\geq\frac{1}{4}-o(1)$ w.h.p.\ in the WP model of closed surfaces. This shows that the optimal spectral gap property is in fact the typical behavior, in the sense of the Weil-Petersson measure, of closed hyperbolic surfaces as $g\to\infty$.

Despite the above advances, the question whether optimal spectral gaps  $\lambda_{1}\geq\frac{1}{4}-o(1)$ are typical for random covers has remained open.
The reason is that strong convergence of uniform random homomorphisms of surface groups to $S_n$ was outside the reach of previous methods. Thus \cite{zbMATH07734891,zbMATH07974824}  had to rely on reductions to the case of free groups: either by working with non-compact surfaces, or by embedding surface groups into an extension of centralizer of a free group which results in a complicated randomization that is far from uniform. The reason we can resolve this question here (Theorem~\ref{thm:main-hyperbolic}) is that the new methods introduced in this paper resolve the corresponding strong convergence problem (Theorem~\ref{thm:strong_convergence}).

We emphasize that all prior results on optimal spectral gaps of surfaces are either not effective, or yield at best logarithmic convergence rates (see, e.g., \cite{arXiv:2305.04584}). This is in sharp contrast 
to the analogous setting of random regular graphs, where the exact Tracy-Widom asymptotics (at scale $\sim n^{-2/3}$) of
the first nontrivial eigenvalue was recently established in a remarkable breakthrough of
Huang--McKenzie--Yau \cite{huang2025ramanujanpropertyedgeuniversality}. Establishing such a  result for random surfaces remains far out of reach using currently known methods. Nonetheless, the
polynomial convergence rate $\sim n^{-c}$ that follows from the present paper and \cite{hidethomasmacera2025} may be viewed as a step in this direction.

\medskip

\emph{Beyond constant curvature.} As one primary interest in hyperbolic surfaces arises from their role as models of quantum chaos, it is desirable to follow up successes in this setting with the corresponding results in variable negative curvature, e.g., \cite{zbMATH02019764, zbMATH05249493, zbMATH07485049}, in order to show the robustness of the relevant phenomenon --- in this case near optimal spectral gaps w.h.p. One benefit of the random cover model is that it applies as-is to variable negative curvature settings. Indeed, as is shown in \cite{hidemoynaud2025}, Theorem \ref{thm:strong_convergence} means that we now have a model of random large genus variable negative curvature surfaces having near optimal spectral gaps with high probability.

\subsection{Paper organization}

The remainder of this paper is organized as follows. In \S\ref{sec:markov}, we state some basic inequalities for polynomials and rational functions that will be used in the sequel.
The proofs of Theorems \ref{thm:main}, \ref{thm:asymp-effective}, and \ref{thm:powers} are subsequently developed in \S\ref{sec:Proof-of-new-criterion-strong-convergence}, \S\ref{sec:EffectivePolynomialApproximation}, and \S\ref{sec:Proper-powers}, respectively.
Finally, we assemble the above ingredients in \S\ref{sec:Main-proof} to complete the proof of Theorem \ref{thm:strong_convergence}.

\subsection{Notation}
\label{sec:notation}

We write $\N=\{1,2,3,\ldots\}$, $\Z_+= \{0,1,2,\ldots\}$, $\Z= \{\ldots,-1,0,1,\ldots\}$, and for $n\in\N$, $[n]=\{1,\ldots,n\}$. The falling Pochhammer symbol and double factorial are
\[
(n)_{m}\eqdf n(n-1)\cdots(n-m+1),\qquad\quad   (2n-1)!!\eqdf (2n-1)(2n-3)\cdots 3\cdot 1,
\]
where by convention $(2n-1)!!=1$ for $n=0$.

A set partition of $[n]$ is a collection of disjoint subsets whose union is $[n]$. We write $\Part([n])$ for the set partitions of $[n]$. 

For $f,g:\N\to\R$, we write $f=O(g)$ if there is a universal constant $C>0$ --- depending on no other parameters --- such that for all $n\in\N,$ $|f(n)|\leq Cg(n)$. If $C$ depends on other parameters this is indicated with a subscript to the big $O$.

We denote by $f^{(m)}$ the $m$th derivative of a univariate function $f$, and denote by $\|f\|_{[a,b]} = \sup_{z\in[a,b]}|f(z)|$ the uniform norm on the interval $[a,b]$.

\subsubsection*{Groups and $C^*$-algebras}

Let $\Lambda$ be a finitely generated group with a finite generating
set $S$. We denote by $|\gamma|$ the word length of $\gamma\in\Lambda$ with 
respect to $S$.

The group algebra $\C[\Lambda]$ has a natural involution
$(\sum_\gamma \alpha_\gamma \gamma)^* = \sum_\gamma \bar\alpha_\gamma \gamma^{-1}$. Thus we can speak
of self-adjoint elements of the group algebra; these become self-adjoint
operators in any unitary representation. For any univariate polynomial
$h$ and $x\in \C[\Lambda]$, we interpret $h(x)$ in the obvious algebraic
sense. We denote by $|x|$ the maximum word length of a group element in 
the support of $x$ when $x$  is in either $\C[\Lambda]$ or the tensor product of $\C[\Lambda]$ with a matrix algebra.  We define 
$$
        \|x\|_{C^*(\Lambda)} = \sup_\pi \|\pi(x)\|
$$
for $x\in \C[\Lambda]$, where the supremum is taken over all unitary 
representations of $\Lambda$.

Let $\lambda$ be the regular representation of $\Lambda$, and extend it linearly to 
$\C[\Lambda]$. The reduced $C^*$-algebra $C^*_{\rm red}(\Lambda)$ of $\Lambda$ is the
norm-closure of $\{\lambda(x):x\in \C[\Lambda]\}$.

For $M\in\mathrm{M}_n(\C)$, let $\tr M$ be the trace and $\ntr 
M=\frac{1}{n}\tr M$ be the normalized trace. We denote by $\tau$ the 
canonical trace on $C^*_{\rm red}(\Lambda)$, i.e., $\tau(a)=\langle\delta_e,
a\delta_e\rangle$, and by
$$
        \|a\|_{L^2(\tau)}=\tau(a^*a)^{1/2}
$$
for $a\in C^*_{\rm red}(\Lambda)$. In particular,
$\|\lambda(x)\|_{L^2(\tau)} = (\sum_\gamma |\alpha_\gamma|^2)^{1/2}$ for
$x=\sum_\gamma \alpha_\gamma \gamma$.

\subsection*{Acknowledgments}

We thank Will Hide, Laura Monk, Bram Petri, Mikael de la Salle, and Joe Thomas for useful discussions. We are grateful to the referee for helpful comments and suggestions on the manuscript.

This work was conducted while D.P.\ and R.v.H.\ were members of the 
Institute for Advanced Study in Princeton, NJ, which is gratefully 
acknowledged for providing a fantastic mathematical environment.

Funding (M.M.): This project has received funding from the European
Research Council (ERC) under the European Union’s Horizon 2020 research
and innovation programme (grant agreement No.\ 949143). This work
was funded in part by a Philip Leverhulme prize awarded by the Leverhulme
Trust.

Funding (D.P.): This work was supported by the European Research Council (ERC) under the European Union’s Horizon 2020 research and innovation programme (grant agreement No.\ 850956), by the Israel Science Foundation, ISF grants 1140/23, as well as the National Science Foundation under grant No.\ DMS-1926686.

Funding (R.v.H.): This work was supported in part by NSF grant DMS-2347954.

\section{Analysis of polynomials and rational functions\label{sec:markov}}

We begin by recalling an elaboration of the Markov brothers inequality that is the basic
principle used in the polynomial method \cite{chen2024new}. The following convenient formulation
may be found in \cite[Lem.\ 4.2]{MageeSalle2}.

\begin{lem}
\label{lem:markov-brothers-interpolated}For every real polynomial $P$ of degree at most $q$ and every $k\in \Z_+$,
\[
\sup_{t\in\left[0,\frac{1}{2q^{2}}\right]}|P^{(k)}(t)|\leq\frac{2^{2k+1}q^{4k}}{(2k-1)!!}\sup_{n\geq q^{2}}\left|P\left(\frac{1}{n}\right)\right|.
\]
\end{lem}
We will require a variant of this principle for rational functions.

\begin{lem}
\label{lem:rational-final}
Let $P,Q$ be real polynomials of degree at most $q\in\N$, and let 
\[
\Phi(t)=\frac{P(t)}{Q(t)}.
\]
Assume that there is a constant $C>1$ so that 
\[
\bigg|P\bigg(\frac{1}{n}\bigg)\bigg|\le C,\qquad\qquad C^{-1}\le Q\bigg(\frac{1}{n}\bigg)\le C
\]
for all $n\in\N$ with $n\ge(Cq)^{C}$. Then there is $C'=C'(C)>1$
so that for all $k\le Cq$ 
\[
\sup_{t\in[0,(C'q)^{-C'}]}|\Phi^{(k)}(t)|\le(C'q)^{C'k}.
\]
\end{lem}

\begin{proof}
Applying Lemma \ref{lem:markov-brothers-interpolated} with $q \leftarrow (Cq)^C$ yields 
\[
\sup_{t\in[0,(2Cq)^{-2C}]}|P^{(k)}(t)|\le2C(2Cq)^{4Ck},\qquad\sup_{t\in[0,(2Cq)^{-2C}]}|Q^{(k)}(t)|\le2C(2Cq)^{4Ck},
\]
for all $k\in\N$. On the other hand, every $t\in[0,(2Cq)^{-2C-1}]$
is within distance at most $(2Cq)^{-4C-2}$ from some point $\frac{1}{n}$
with $n\ge(Cq)^{C}$. Thus 
\[
\inf_{t\in[0,(2Cq)^{-2C-1}]}Q(t)\ge C^{-1}-(2Cq)^{-4C-2}\sup_{t\in[0,(2Cq)^{-2C-1}]}|Q^{(1)}(t)|\ge(2C)^{-1}.
\]
The conclusion now follows from the chain rule. Indeed, we can first estimate
\[
    \sup_{t\in[0,(2Cq)^{-2C-1}]}\left|\left(\frac{1}{Q}\right)^{(k)}(t)\right| \le (2C)^{2k+1} (2Cq)^{4Ck} k!
\]
as in the last part of \cite[proof of Lemma 7.2]{MageeSalle2}. Thus
\[
    |\Phi^{(k)}(t)|
    =
    \left|\sum_{i=0}^k \binom{k}{i} P^{(k-i)}(t)
    \left(\frac{1}{Q}\right)^{(i)}(t)\right| \le 2^k (2C)^{2k+2} (2Cq)^{4Ck} k!
\]
for $t\in[0,(2Cq)^{-2C-1}]$.
The conclusion follows as $k!\le (Cq)^k$ for $k\le Cq$.
\end{proof}

We finally recall a simple fact about trigonometric polynomials.

\begin{lem}
\label{lem:fourier} Let $f(\theta)=\sum_{k=0}^{q}a_{k}\cos(k\theta)$ and $m\in\N$. Then 
\[
\sum_{k=1}^{q}k^{m-1}|a_{k}|\le4\|f^{(m)}\|_{[0,2\pi]}.
\]
\end{lem}

\begin{proof}
As $(k^m a_k)_{0\le k\le q}$ are the Fourier coefficients of $f^{(m)}$, we have
\[
    \sum_{k=1}^{q}k^{m-1}|a_{k}|\le
    \left(\sum_{k=1}^{q} \frac{1}{k^2}\right)^{1/2}
    \left(
    \sum_{k=1}^{q}|k^{m}a_{k}|^2
    \right)^{1/2}
    \le
    \sqrt{\frac{\pi^2}{6}} \| f^{(m)}\|_{L^2[0,2\pi]}
\]
by Cauchy-Schwarz and Parseval. We conclude using $\|f\|_{L^2[0,2\pi]}\le\sqrt{2\pi}\|f\|_{[0,2\pi]}$.
\end{proof}

\section{The generalized polynomial method\label{sec:Proof-of-new-criterion-strong-convergence}}

The aim of this section is to prove Theorem \ref{thm:main}. The reader who is new to the polynomial method is encouraged to review \cite[\S 2]{chen2024new} for an introduction, before proceeding to the details of the proof.

\subsection{Proof of Theorem \ref{thm:main}}

In this section, we always assume that Assumptions \ref{ass:cond1} and \ref{ass:cond2} hold. We denote by $C$ the constant in Assumption \ref{ass:cond1}, and by $r=|S|$ the size of the generating set of $\Lambda$.

The key step in the proof is to establish the following master inequality.

\begin{lem}
\label{lem:key} 
Fix a self-adjoint $x\in\C[\Lambda]$ and
let $K=\|x\|_{C^{*}(\Lambda)}$. Then 
\[
|u_{1}(h(x))|\le r(6Cq|x|)^{6C}\|h\|_{[-K,K]}
\]
and 
\[
\bigg|\E\big[{\ntr h(\pi_{n}(x))}\big]-\tau\big(h(\lambda(x))\big)-\frac{1}{n}u_{1}(h(x))\bigg|\le\frac{2r^{2}(6Cq|x|)^{12C}}{n^{2}}\|h\|_{[-K,K]}
\]
for all $n\ge1$ and every real polynomial $h$ of degree at most $q$. 
\end{lem}

\begin{proof}
We assume without loss of generality that $|x|\ge1$ and $q\ge1$ (otherwise $h(x)$ is a multiple of the identity, so both sides of the inequalities vanish). Let us write 
\[
h(x)=\sum_{\gamma\in\Lambda}\alpha_{\gamma}\gamma.
\]
As $|h(x)|\le q|x|$, there are at most $(2r+1)^{q|x|}\le(3r)^{q|x|}$ nonzero coefficients $\alpha_{g}$. Applying Assumption \ref{ass:cond1} separately to each element of the support of $h(x)$ yields 
\[
\Bigg|\E\big[{\ntr h(\pi_{n}(x))}\big]-f_{h}\bigg(\frac{1}{n}\bigg)\Bigg|\le\frac{(3Cq|x|)^{3Cq|x|}}{n^{q|x|+2}}\sum_{\gamma\in\Lambda}|\alpha_{\gamma}|
\]
for all $n\ge C(3q|x|)^{C}$, where we define the polynomial 
\[
f_{h}(t)=\sum_{k=0}^{q|x|+1}u_{k}(h(x))\,t^{k}
\]
and we used that $q|x|+2\le3q|x|$. As 
\begin{equation}
\label{eq:ignorecancellations}
\sum_{\gamma}|\alpha_{\gamma}|\le(3r)^{q|x|/2}\|h(\lambda(x))\|_{L^{2}(\tau)}\le(3r)^{q|x|/2}\|h\|_{[-K,K]}
\end{equation}
by Cauchy-Schwarz and as $\|\lambda(x)\|\le K$, we readily estimate
\begin{equation}
\Bigg|\E\big[{\ntr h(\pi_{n}(x))}\big]-f_{h}\bigg(\frac{1}{n}\bigg)\Bigg|\le\frac{1}{n^{2}}\|h\|_{[-K,K]}\label{eq:exp}
\end{equation}
for $n\ge(3r)^{1/2}(3Cq|x|)^{3C}$. Now note that (\ref{eq:exp}) implies the \emph{a priori} bound 
\[
\bigg|f_{h}\bigg(\frac{1}{n}\bigg)\bigg|\le\big|\E\big[{\ntr h(\pi_{n}(x))}\big]\big|+\frac{1}{n^{2}}\|h\|_{[-K,K]}\le2\|h\|_{[-K,K]}.
\]
Moreover, $u_{k}(h(x))$ are real as $h(x)$ is self-adjoint, so $f_{h}$ is a real polynomial. We can therefore apply Lemma \ref{lem:markov-brothers-interpolated} with $q^2\leftarrow (3r)^{1/2}(3Cq|x|)^{3C}$ to estimate 
\begin{align*}
 & \|f_{h}'\|_{[0,r^{-1/2}(6Cq|x|)^{-3C}]}\le r(6Cq|x|)^{6C}\|h\|_{[-K,K]},\\
 & \|f_{h}''\|_{[0,r^{-1/2}(6Cq|x|)^{-3C}]}\le r^{2}(6Cq|x|)^{12C}\|h\|_{[-K,K]}.
\end{align*}
We can now readily conclude the proof. Indeed, as $u_{1}(h(x))=f_{h}'(0)$,
the first part of the statement is immediate. On the other hand, Taylor expanding $f_h$ yields
\[
\bigg|f_{h}\bigg(\frac{1}{n}\bigg)-u_{0}(h(x))-\frac{1}{n}u_{1}(h(x))\bigg|\le\frac{1}{n^{2}}\|f_{h}''\|_{[0,\frac{1}{n}]}\le\frac{r^{2}(6Cq|x|)^{12C}}{n^{2}}\|h\|_{[-K,K]}
\]
for $n\ge r^{1/2}(6Cq|x|)^{3C}$, so the second part follows in this
case using (\ref{eq:exp}) and that $u_{0}(x)=\tau(\lambda(x))$ by
Assumption \ref{ass:cond1}. Finally, when $n<r^{1/2}(6Cq|x|)^{3C}$,
we have 
\begin{align*}
\bigg|\E\big[{\ntr h(\pi_{n}(x))}\big]-\tau\big(h(\lambda(x))\big)-\frac{1}{n}u_{1}(h(x))\bigg| & \le2\|h\|_{[-K,K]}+\frac{1}{n}|u_{1}(h(x))|\\
 & \le\frac{2r^{2}(6Cq|x|)^{12C}}{n^{2}}\|h\|_{[-K,K]}
\end{align*}
using the triangle inequality, the bound on $|u_{1}(h(x))|$, and
$1\le\frac{r^{1/2}(6Cq|x|)^{3C}}{n}$. 
\end{proof}
With the statement of Lemma \ref{lem:key} in hand, the remainder
of the argument is now precisely as in \cite{chen2024new}. We first
extend Lemma \ref{lem:key} to smooth $h$.
\begin{cor}
\label{cor:smaster} Fix a self-adjoint $x\in\C[\Lambda]$
and let $K=\|x\|_{C^{*}(\Lambda)}$. Then $\nu(h)\eqdf u_{1}(h(x))$ extends
to a compactly supported distribution, and 
\[
\bigg|\E\big[{\ntr h(\pi_{n}(x))}\big]-\tau\big(h(\lambda(x))\big)-\frac{1}{n}\nu(h)\bigg|\le\frac{8r^{2}(6C|x|)^{12C}}{n^{2}}\|f^{(m)}\|_{[0,2\pi]}
\]
for all $n\ge1$ and $h\in C^{\infty}(\R)$. Here $m=1+\lceil12C\rceil$
and $f(\theta)=h(K\cos\theta)$. 
\end{cor}

\begin{proof}
That $\nu$ extends to a compactly supported distribution follows
directly from the first part of Lemma \ref{lem:key} and \cite[Lem.\ 4.7]{chen2024new}.

Now assume first that $h$ is a polynomial of degree $q$. Then we
can uniquely write 
\[
h(t)=\sum_{k=0}^{q}a_{k}T_{k}\bigg(\frac{t}{K}\bigg),
\]
where $T_{k}$ is the Chebyshev polynomial defined by $T_{k}(\cos(\theta))=\cos(k\theta)$.
Applying the second part of Lemma \ref{lem:key} separately to each
Chebyshev polynomial yields 
\[
\bigg|\E\big[{\ntr h(\pi_{n}(x))}\big]-\tau\big(h(\lambda(x))\big)-\frac{1}{n}u_{1}(h(x))\bigg|\le\frac{2r^{2}(6C|x|)^{12C}}{n^{2}}\sum_{k=1}^{q}k^{12C}|a_{k}|.
\]
The conclusion now follows from Lemma \ref{lem:fourier} in the case
that $h$ is a polynomial. Since the resulting bound does not depend
on the degree of $h$ and as polynomials are dense in $C^{\infty}(\R)$,
the conclusion extends by continuity to any $h\in C^{\infty}(\R)$. 
\end{proof}
We can now state a quantitative version of Theorem \ref{thm:main}.
\begin{thm}
\label{thm:quantmain} Fix a self-adjoint $x\in\C[\Lambda]$
and let $K=\|x\|_{C^{*}(\Lambda)}$. Then 
\[
\P\big[\|\pi_{n}(x)\|\ge\|\lambda(x)\|+\varepsilon\big]\le\frac{r^{2}}{n}\bigg(\frac{c|x|K}{\varepsilon}\bigg)^{c}
\]
for every $n\ge1$ and $\varepsilon>0$. Here $c$ is a constant that
depends only on $C$. 
\end{thm}

\begin{proof}
By \cite[Lem.\ 4.10]{chen2024new}, there exists $h\in C^{\infty}(\R)$
with values in $[0,1]$ so that 
\begin{enumerate}
\item $h(z)=0$ for $|z|\le\|\lambda(x)\|+\frac{\varepsilon}{2}$ and $h(z)=1$
for $|z|\ge\|\lambda(x)\|+\varepsilon$; 
\item $\|f^{(m)}\|_{[0,2\pi]}\le(\frac{cK}{\varepsilon})^{c}$ for $f(\theta)=h(K\cos\theta)$. 
\end{enumerate}
Moreover, Assumption \ref{ass:cond2} and \cite[Lem.\ 4.9]{chen2024new}
yield $\supp\nu\subseteq[-\|\lambda(x)\|,\|\lambda(x)\|]$, where
$\nu$ is defined in Corollary \ref{cor:smaster}. Thus $\tau(h(\lambda(x)))=\nu(h)=0$,
so that 
\[
\P\big[\|\pi_{n}(x)\|\ge\|\lambda(x)\|+\varepsilon\big]\le\E\big[{\tr h(\pi_{n}(x))}\big]\le\frac{8r^{2}(6C|x|)^{12C}}{n}\bigg(\frac{cK}{\varepsilon}\bigg)^{c}
\]
using Corollary \ref{cor:smaster}. This concludes the proof. 
\end{proof}
The statement of Theorem \ref{thm:main} now follows immediately from
Theorem \ref{thm:quantmain} (if $x$ is not self-adjoint, we simply
apply Theorem \ref{thm:quantmain} to $x^{*}x$). However, Theorem
\ref{thm:quantmain} provides much stronger information, since it
even yields a polynomial rate: for all $b<1/c$
\[
\|\pi_{n}(x)\|\le\|\lambda(x)\|+O(n^{-b}) \quad\text{with probability}\quad 1-o(1).
\]

\subsection{Remarks and complements}

\hspace*{\parindent}1. Assumption \ref{ass:cond1} ensures that $\E[\ntr h(\pi_n(x))]$ can be approximated by a polynomial of $\frac{1}{n}$. At first sight, however, it is somewhat surprising that this suffices for the polynomial method. As is explained in \cite[\S 2.3]{chen2024new}, the key feature of this method is that it captures certain delicate cancellations between the coefficients of $h(x)=\sum_\gamma \alpha_\gamma \gamma$. On the other hand, Assumption \ref{ass:cond1} only applies to each $\gamma$ separately and is therefore unable to preserve cancellations.

The reason this does not matter is that the price for disregarding cancellations is exponential in $q|x|$ (see \eqref{eq:ignorecancellations}), and can therefore be absorbed 
by the error term in Assumption \ref{ass:cond1} by expanding to order $\sim q|x|$. We can therefore first approximate the spectral statistics by high 
degree polynomials in the most naive manner, and then proceed to 
capture the cancellations in the latter using the polynomial method. This is crucial for applications, since we do not know how to establish \emph{a priori} bounds as in Assumption \ref{ass:cond1} that capture cancellations.

\medskip

2. While we developed the argument above in the group setting, the method
is very robust and applies to many other situations. The only part
of the argument that relied on the group structure is \eqref{eq:ignorecancellations}.
In other situations,
one can achieve a similar conclusion by applying a classical result
of V.\ Markov which states that 
\[
\sum_{k=0}^{q}|a_{k}|\le e^{q/K}\|h\|_{[-K,K]}
\]
for every real polynomial $h(t)=\sum_{k=0}^{q}a_{k}t^{k}$ and $K>0$ \cite[\S  2.6, Eq.\ (9)]{zbMATH03191390},
at the expense of a less explicit dependence of the constants on the
choice of $x$.

\medskip

3. While we have proved strong convergence of $\pi_n(x)$ for $x\in\C[\Lambda]$ for simplicity of exposition, the proof is readily modified to establish strong
convergence with matrix coefficients, that is, of $[\mathrm{id}\otimes\pi_n](x)$ for $x\in \mathrm{M}_d(\C)\otimes\C[\Lambda]$.
In this setting, Assumption~\ref{ass:cond2} must be replaced by the assumption that $u_1$ is ``completely tempered'' \cite{MageeSalle2}.

\begin{assumption}
\label{ass:cond2bis}
For every $d\in\N$ and self-adjoint $x\in\mathrm{M}_d(\C)\otimes\C[\Lambda]$, we have
\[
        \limsup_{p\to\infty} |[\mathrm{tr}\otimes u_1](x^p)|^{1/p} \le
        \|[\mathrm{id}\otimes\lambda](x)\|.
\]
\end{assumption}

\begin{thm}
\label{thm:quantmtx}
Suppose that Assumptions \ref{ass:cond1} and \ref{ass:cond2bis} hold. Then
for any self-adjoint $x\in \mathrm{M}_d(\C)\otimes\C[\Lambda]$ and for every $n\ge1$ and $\varepsilon>0$,
we have
\[
\P\big[\|[\mathrm{id}\otimes\pi_{n}](x)\|\ge(1+\varepsilon)\|[\mathrm{id}\otimes\lambda](x)\|\big]\le\frac{cd}{n\varepsilon^b}.
\]
Here $b$ is a constant that depends only on $C$, and $c$ depends on
$C$, $r$, and $|x|$.
\end{thm}

\begin{proof}
We describe the requisite changes to the proof of Theorem \ref{thm:quantmain}.
In the proof of Lemma \ref{lem:key}, we can still write
$h(x)=\sum_{\gamma\in\Lambda}\alpha_{\gamma}\gamma$ where now
$\alpha_\gamma \in \mathrm{M}_d(\C)$ are matrix coefficients. Instead of
\eqref{eq:ignorecancellations}, we may now estimate
\begin{equation}
\label{eq:passtoreducednorm}
    \sum_{\gamma}\|\alpha_{\gamma}\|\le(3r)^{q|x|} 
    \max_\gamma \|\alpha_\gamma\| \le  
    (3r)^{q|x|} \|h([\mathrm{id}\otimes\lambda](x))\|
    \le(3r)^{q|x|}\|h\|_{[-K,K]},
\end{equation}
where we used \cite[eq.~(9.7.2)]{pisier2003operatorspacebook} in the second inequality.
Repeating the rest of the proof of Theorem \ref{thm:quantmain} with only cosmetic changes
yields
\[
\P\big[\|[\mathrm{id}\otimes\pi_{n}](x)\|\ge\|[\mathrm{id}\otimes\lambda](x)\|+\delta\big]\le\frac{d}{n}\bigg(\frac{c'K}{\delta}\bigg)^{b},
\]
where $b$ depends only on $C$ and $c'$ depends on $C,r,|x|$ (note that we multiplied in the proof of Theorem \ref{thm:quantmain} by $dn$, rather than $n$, to pass from the normalized to the unnormalized trace).
We finally note that $K=\|x\|_{\mathrm{M}_d(\C)\otimes C^*(\Lambda)}\le (3r)^{|x|}\|[\mathrm{id}\otimes\lambda](x)\|$
by the same argument as in \eqref{eq:passtoreducednorm}, so the
conclusion follows with $\delta\leftarrow \varepsilon\|[\mathrm{id}\otimes\lambda](x)\|$.
\end{proof}

Theorem \ref{thm:quantmtx} shows that a strong convergence upper bound with a polynomial rate holds when $|x|$ is bounded and we admit matrix coefficients of dimension $d \le n^{1-\beta}$. This is used in \cite{hidethomasmacera2025} to obtain a polynomial rate in Theorem \ref{thm:main-hyperbolic}.

\section{Effective approximation by polynomials}
\label{sec:EffectivePolynomialApproximation}

The aim of this section is to prove Theorem \ref{thm:asymp-effective}.

\emph{Throughout this section, $C>0$ will denote a constant that may depend on the genus $g$,
which can change from line to line but will only increase finitely many times.}

\subsection{Background}

\subsubsection*{Representation theory of the symmetric group}

A (number\footnote{As opposed to a set partition.}) partition of
$n$ is a non-increasing sequence of non-negative integers\\
$\lambda=(\lambda_{1},\lambda_{2},\ldots,\lambda_{\ell(\lambda)})$
with 
\[
\sum_{i=1}^{\ell(\lambda)}\lambda_{i}=n.
\]
In this case we write $\lambda\vdash n$ and say $\lambda$ has size
$n$, $|\lambda|=n$. We think of partitions and Young diagrams (YDs)
interchangeably in the sequel. If $\lambda$ is a YD we write $\lambda^{\vee}$
for the conjugate YD obtained by interchanging rows and columns. We
write $b_{\lambda}$ for the number of boxes outside the first row
of $\lambda$ and $b_{\lambda}^{\vee}\eqdf b_{\lambda^{\vee}}$ the
number of boxes outside the first column.

The equivalence classes of irreducible representations (irreps) of
$S_{n}$ are parameterized (using Young symmetrizers, for example)
by the YDs $\lambda\vdash n$. Given $\lambda\vdash n,$ we write
$\chi_{\lambda}$ for the character of the associated irrep and $d_{\lambda}=\chi_{\lambda}(1)$
for its dimension. We write $V^{\lambda}$ for the corresponding $S_{n}$-module.
As $\chi_{\lambda^{\vee}}=\mathrm{sign}\cdot\chi_{\lambda}$ we have
$d_{\lambda}=d_{\lambda^{\vee}}$.

If one YD $\mu$ is contained in another $\lambda$, we say $\mu\subset\lambda$,
and if they differ by $k$ boxes, we write $\mu\subset_{k}\lambda$.
For YDs $\mu\subset_{k}\lambda$, a (standard) tableau of shape $\lambda/\mu$
is a filling of the $k$ boxes of $\lambda$ not in $\mu$ by the
numbers in $[k]$ such that the numbers in each row (resp.~column) are increasing from left to right (resp.~top to bottom). We write $\Tab(\lambda/\mu)$ for the collection of such tableaux. It is obvious that when $\mu\subset_{k}\lambda$,
\[
|\Tab(\lambda/\mu)|\leq k!.
\]
We use the fact that every irrep of $S_{n}$ can be realized over
the reals and therefore as orthogonal matrices (so no complex conjugates
appear in our calculations).

\subsubsection*{The Witten zeta function}

The Witten zeta function of $S_{n}$ is defined by 
\[
\zeta(s;S_{n})\eqdf\sum_{\lambda\vdash n}\frac{1}{d_{\lambda}^{s}}.
\]
This appears as a normalization factor in the uniform measure on $\Hom(\Gamma_{g},S_{n})$
and hence is important in the sequel. By a result of Liebeck and Shalev
\cite[Prop.~2.5]{LiebeckShalev} (independently, Gamburd \cite[Prop.~4.2]{GamburdBelyi}),
for real $s>0$,
\begin{equation}
\zeta(s;S_{n})=2+O_{s}\left(n^{-s}\right)\label{eq:zeta-limit}
\end{equation}
 as $n\to\infty$. We also have the following effective tail bound  \cite[Prop.~4.6]{MageeNaudPuder}.
 
\begin{prop}
\label{prop:mageenaudpuder}
For all $s>0$, there is $\kappa=\kappa(s)>1$ such that 
\begin{equation}
\sum_{\substack{\lambda\vdash n\\
b_{\lambda},b_{\lambda}^{\vee}\geq b
}
}\frac{1}{d_{\lambda}^{s}}\leq\left(\frac{\kappa b^{2s}}{(n-b^{2})^{s}}\right)^{b}
\label{eq:zeta_tail}
\end{equation}
for all $b,n\in\N$ with $b^{2}\leq\frac{n}{3}$.
\end{prop}

\subsubsection*{Expansion of $\protect\E_{g,n}[\protect\fix_{\gamma}]$ by expected
numbers of embeddings}

In the following, $\phi$ is a uniformly random element of $\Hom(\Gamma_{g},S_{n})$.
We now recap on the language of \cite{MPasympcover}. This is just
a formalism for expressing the facts:
\begin{itemize}
\item $\E_{g,n}[\fix_{\gamma}]$ is $n$ times the probability that $\phi(\gamma)$
fixes a given point, say $1\in[n]$, by invariance under $S_{n}$
conjugation.
\item The event $\mathcal{E}_{\gamma}$ that $\phi(\gamma)$ fixes $1$
can be partitioned into sub-events in a variety of different ways. 
\end{itemize}
In \cite{MPasympcover} these sub-events are chosen more carefully
than we have to here. 

We use the standard set of generators $a_1,\ldots,b_g$ of $\Gamma_{g}$ from \eqref{eq:presentation of Gamma}.  The partitioning is described by reference to the Schreier graph 
\[
X_{\phi}\eqdf\mathrm{Schreier}\left(\{\phi(a_{1}),\ldots,\phi(b_{g})\},[n]\right)
\]
defined by the action of $\phi(a_1),\ldots,\phi(b_g)$ on $[n]$,
i.e., there is a directed edge from $i$ to $j$ for each $f\in\{a_1,\ldots,b_g\}$ such that 
\[
\phi(f)[i]=j.
\]
Every directed edge is colored by the generator from which it arises.\footnote{The definition of $X_{\phi}$ given here is the 1-skeleton of
the $X_{\phi}$ of \cite{MPasympcover}, which includes $n$ additional
2-cells to form a closed surface.}

Given $\gamma\in\Gamma_{g}$, let $\gamma'$ be a conjugate of $\gamma$
whose shortest representing word in the generators is cyclically reduced.
We have $|\gamma'|\leq|\gamma|$ and
\[
\E_{g,n}[\fix_{\gamma}]=\E_{g,n}[\fix_{\gamma'}],
\]
so without loss of generality in the following suppose $\gamma$ is
cyclically reduced.

Let $F_{2g}$ denote the free group on $a_1,\ldots,b_g$. Let
$C_{\gamma}$ denote a circle subdivided into $|\gamma|$ intervals
(viewed as an undirected graph). Direct and label these intervals
by $a_1,\ldots,b_g$ according to an expression of $\gamma$
of length $|\gamma|$ in these generators --- the direction indicates
whether a generator or its inverse appears at a given location. 

Let $\mathcal{R}$ denote the collection of surjective labeled-graph
morphisms
\[
r:C_{\gamma}\twoheadrightarrow W_{r}
\]
such that:
\begin{enumerate}
\item $W_{r}$ is \emph{folded} in the sense that every vertex has at most one incoming $f$-labeled half-edge and at most one outgoing $f$-labeled half-edge, for each $f\in\{a_1,\ldots,b_g\}$. (Note that the assumption $\gamma$
is cyclically reduced means $C_{\gamma}$ itself is folded.)
\item Every path in $W_{r}$ spelling out an element of $F_{2g}$ that is
in the kernel of $F_{2g}\to\Gamma_{g}$ is closed.
\end{enumerate}
In the language of \cite{MageePuderCore,MPasympcover} these $W_{r}$
are (very special cases of) \emph{tiled surfaces}.\footnote{In (ibid.) it is a canonical thickening of the graph that is thought
of as the surface.}

Any homomorphism of folded labeled graphs $C_{\gamma}\to X_{\phi}$
factors uniquely as a surjective morphism followed by an injective
one, namely
\[
C_{\gamma}\stackrel{r}{\twoheadrightarrow}W_{r}\hookrightarrow X_{\phi},
\]
for unique $r\in\mathcal{R}$. In the language of \cite{MPasympcover},
this factoring property means $\mathcal{R}$ is a \emph{resolution}
of $C_{\gamma}$. 
\begin{rem}[Comparison to \cite{MPasympcover}]
Although the resolution we define here might be the most natural,
the methods of \cite{MPasympcover} leading to \cite[Thm.~1.2]{MPasympcover}
--- not proved here --- only work on $W_{r}$ that have, very roughly
speaking, geodesic boundary. So in \cite{MPasympcover} great care
is taken over finding a finer resolution whose boundary components
have desired properties. For our purposes, the simplest resolution suffices.
\end{rem}

We have 
\begin{equation}
|\mathcal{R}|\leq|\gamma|!\label{eq:resolution-bound}
\end{equation}
since each element is defined by the partition of the vertices of
$C_{\gamma}$ given by the fibers of the surjective immersion, and
there are $|\gamma|$ vertices. By \cite[Lem.\ 2.7 and 2.9]{MPasympcover}
\begin{equation}
\E_{g,n}[\fix_{\gamma}]=\sum_{r\in\mathcal{R}}\E_{n}^{\mathrm{\emb}}(W_{r}), \label{eq:resolution-expansion}
\end{equation}
where 
\[
\E_{n}^{\emb}(W_{r})=\E_{\phi\in\X_{g,n}}\#\{\,\text{injective morphisms \ensuremath{W_{r}\to X_{\phi}\,\}.}}
\]

In the sequel, for $W_{r}$ as above we write $\v(W_{r}),\e(W_{r})$
for the number of vertices and edges of $W_{r}$. For $f\in\{a_1,\ldots,b_g\}$
we write $\e_{f}(W_{r})$ for the number of $f$-labeled edges of
$W_{r}$. With things as before, clearly 
\[
\v(W_{r}),\e(W_{r}),\e_{f}(W_{r})\leq|\gamma|.
\]

\subsubsection*{Expected number of embeddings, prior results}

\emph{In the remainder of \S\ref{sec:EffectivePolynomialApproximation}, we restrict to $g=2$ for simplicity of exposition.
We now use $a,b,c,d$ for $a_1,b_1,a_2,b_2$. The only things that depend on $g$ in a non-obvious way are constants and we indicate how these depend on $g$ throughout. All integrals in the rest of this section are taken with respect to uniform probability measures.}

\medskip

Let $Y$ denote some fixed $W_{r}$ from our resolution $\mathcal{R}$.
Let $\v=\v(Y)$, $\e_{f}=\e_{f}(Y)$. In the following we label the
vertices of $Y$ injectively by $1,\ldots,\v$.

For $g_{f}\in S_{n}$, $f\in\{a,b,c,d\}$, we say $g_{f}$ \emph{obey}
$\ensuremath{Y}$ if the fixed labeling of the vertices of $Y$ induces
an embedding
\[
Y\hookrightarrow\mathrm{Schreier}\left(\{g_{f}\},[n]\right),
\]
or in simple terms, if $Y$ has an $f$-colored directed edge from
a vertex labeled $i$ to vertex labeled $j$, then $g_{f}(i)=j$.

We have by \cite[(2.1), (5.6), and proof of Prop.~5.1]{MPasympcover}, since $\v(Y)\leq q$, for $n\geq q$
\[
\E_{n}^{\emb}(Y)=\frac{1}{\zeta(2;S_{n})}\frac{(n)_{\v}}{\prod_{f\in a,b,c,d}(n)_{\e_{f}}}\sum_{\lambda\vdash n}d_{\lambda}\Theta_{\lambda}\left(Y\right),
\]
where
\begin{equation}
\Theta_{\lambda}\left(Y\right)\eqdf\left(\prod_{f\in a,b,c,d}(n)_{\e_{f}}\right)\int_{g_{f}\in S_{n}}\mathbf{1}\{g_{f}\text{ obey \ensuremath{Y}\}}\chi_{\lambda}\left([g_{a},g_{b}][g_{c},g_{d}]\right)\label{eq:theta-first}
\end{equation}
(there is a change of normalization since we pass between integrals
over different groups: in the notation of \emph{(ibid.)} $|G_{f}|=(n-\e_{f})!$).
\begin{rem}
In \cite{MPasympcover} the definition of $\Theta_{\lambda}\left(Y\right)$
depends on an auxiliary labeling $\J_{n}:Y^{(0)}\to\{n-\v+1,\ldots,n\}$.
This is not needed here, as the definition of $\Theta_{\lambda}\left(Y\right)$
in (\ref{eq:theta-first}) is the same for any labeling 
by invariance of the uniform measure on $S_n$ by conjugation.
\end{rem}

From \cite[Prop.~5.8]{MPasympcover}
\begin{align}
\Theta_{\lambda}(Y)= & \sum_{\nu\subset_{\v}\lambda}d_{\nu}\sum_{\substack{\nu\subset\mu_{f}\subset_{\e_{f}}\lambda\\
f\in\{a,b,c,d\}
}
}\frac{1}{d_{\mu_{a}}d_{\mu_{b}}d_{\mu_{c}}d_{\mu_{d}}}\Upsilon_{n}\left(\left\{ \sigma_{f}^{\pm},\tau_{f}^{\pm}\right\} ,\nu,\left\{ \mu_{f}\right\} ,\lambda\right),\label{eq:theta-expression-1}
\end{align}
where 
\begin{align}
\Upsilon_{n}\left(\left\{ \sigma_{f}^{\pm},\tau_{f}^{\pm}\right\} ,\nu,\left\{ \mu_{f}\right\} ,\lambda\right)\eqdf\sum_{\begin{gathered}r_{f}^{+},r_{f}^{-}\in\Tab\left(\mu_{f}/\nu\right)\\
s_{f},t_{f}\in\Tab\left(\lambda/\mu_{f}\right)
\end{gathered}
}\M\left(\left\{ \sigma_{f}^{\pm},\tau_{f}^{\pm},r_{f}^{\pm},s_{f},t_{f}\right\} \right)\label{eq:Upsilon-def}
\end{align}
and $\M\left(\left\{ \sigma_{f}^{\pm},\tau_{f}^{\pm},r_{f}^{\pm},s_{f},t_{f}\right\} \right)$
is a product of matrix coefficients of unitary operators on unit vectors
as in \cite[eq.~(5.15)]{MPasympcover}. All we need here about this
product is that 
\begin{equation}
\left|\M\left(\left\{ \sigma_{f}^{\pm},\tau_{f}^{\pm},r_{f}^{\pm},s_{f},t_{f}\right\} \right)\right|\leq1.\label{eq:Mbound}
\end{equation}

\subsection{Tail estimate for $\protect\E_{n}^{\mathrm{\protect\emb}}(W_{r})$}
\label{sec:taileffpoly}

In this section we prove:
\begin{prop}
\label{prop:tail-cut}There is $C>0$ such that for any $q\in \N$, if $\v(Y)\leq q$
and $n\geq28q^{2}$ then
\begin{align*}
\E_{n}^{\emb}(Y) & =\frac{1}{\zeta(2;S_{n})}\frac{(n)_{\v}}{\prod_{f\in a,b,c,d}(n)_{\e_{f}}}\sum_{\substack{\lambda\vdash n\\
\min(b_{\lambda},b_{\lambda}^{\vee})\leq4q
}
}d_{\lambda}\Theta_{\lambda}\left(Y\right)\\
 & \,\,+O\left((Cq)^{Cq}n^{-q}\right).
\end{align*}
\end{prop}

\begin{proof}
Using $\left|\Tab\left(\mu_{f}/\nu\right)\right|,\left|\Tab\left(\lambda/\mu_{f}\right)\right|\leq\v!\leq q!$
together with (\ref{eq:Mbound}) in (\ref{eq:theta-expression-1})
we obtain
\[
\Theta_{\lambda}(Y)\leq(Cq)^{Cq}\sum_{\nu\subset_{\v}\lambda}\frac{1}{d_{\nu}^{3}}.
\]
Note if $b_{\lambda}>4q$ then $b_{\nu}\geq b_{\lambda}-\v>3q$. We
have then 
\[
\sum_{\substack{\lambda\vdash n\\
b_{\lambda},b_{\lambda}^{\vee}>4q
}
}d_{\lambda}\Theta_{\lambda}\left(Y\right)\leq(Cq)^{Cq}\sum_{\lambda\vdash n}\sum_{\substack{\nu\subset_{\v}\lambda\\
b_{\nu},b_{\nu}^{\vee}>3q
}
}\frac{d_{\lambda}}{d_{\nu}^{3}}.
\]
As in \cite[proof of Lemma 5.23]{MPasympcover} the above is 
\begin{equation}
\leq(n)_{\v}(Cq)^{Cq}\sum_{\substack{\nu\vdash n-\v\\
b_{\nu},b_{\nu}^{\vee}>3q
}
}\frac{1}{d_{\nu}^{2}}.\label{eq:one-on-dnu-suqred}
\end{equation}
Hence by (\ref{eq:zeta_tail}) 
\begin{equation}
\sum_{\substack{\nu\vdash n-\v\\
b_{\nu},b_{\nu}^{\vee}>3q
}
}\frac{1}{d_{\nu}^{2}}\leq C\left(\frac{(Cq)^{4}}{(n-\v-9q^{2})^{2}}\right)^{3q}\leq C(Cq)^{Cq}n^{-6q}\label{eq:on-on-dnu-squared-bound}
\end{equation}
in e.g.~$n\geq28q^{2}$\@. (Note that to use (\ref{eq:zeta_tail})
for $b=3q$ as above we need $\v+3(3q)^{2}\leq n$.) Combining all
previous arguments gives
\begin{align*}
\E_{n}^{\emb}(Y) & =\frac{1}{\zeta(2;S_{n})}\frac{(n)_{\v}}{\prod_{f\in a,b,c,d}(n)_{\e_{f}}}\sum_{\substack{\lambda\vdash n\\
\min(b_{\lambda},b_{\lambda}^{\vee})\leq4q
}
}d_{\lambda}\Theta_{\lambda}\left(Y\right)\\
 & +O\left(\frac{1}{\zeta(2;S_{n})}\frac{(n)_{\v}^{2}}{\prod_{f\in a,b,c,d}(n)_{\e_{f}}}(Cq)^{Cq}n^{-6q}\right).
\end{align*}
The denominator Pochhammer symbols in the error are $\geq1$ and numerators
$\leq n^{2q}$. Hence (also using (\ref{eq:zeta-limit})) the whole
error is on the order of
\[
(Cq)^{Cq}n^{-q}
\]
as required.

For general $g\ge2$, the number of tableaux being summed over in
the combination of (\ref{eq:theta-expression-1}) and (\ref{eq:Upsilon-def})
is $8g$, hence $C$ depends linearly on $g$, the summand of (\ref{eq:one-on-dnu-suqred})
is $\frac{1}{d_{\nu}^{2g-2}}\leq\frac{1}{d_{\nu}^{2}}$ and the estimate
(\ref{eq:on-on-dnu-squared-bound}) can be kept the same. The rest
of the proof is the same, changing the product $\prod_{f\in a,b,c,d}(n)_{\e_{f}}$
to a product over $2g$ generators and $\zeta(2;S_{n})$ to $\zeta(2g-2;S_{n}$).
\end{proof}

\subsection{Rationality of contributions from fixed representations}

Given a YD $\lambda\vdash b$, for $n\geq2b$ we write $\lambda^{+}(n)$
for the YD of size $n$ with $\lambda$ outside the first row. Our
aim here is to prove the following rationality result.
\begin{prop}
\label{prop:rationality}Let $b\in\Z_+$, $n\geq\v(Y)+9b$,
and $\lambda\vdash b$. As a function of $n$, $\Theta_{\lambda^{+}(n)}(Y)$
agrees with a rational function of $n$ with coefficients in $\Q$
whose denominator can be taken to be
\[
(n)_{\v+2b}^{5}.
\]
\end{prop}

The proof of Proposition \ref{prop:rationality} uses a
result of Cassidy \cite{cassidy2023projectionformulasrefinementschurweyljones}.
Let $\{e_{i}\}_{i=1}^{n}$ denote the standard orthonormal basis of
$\C^{n}$ with its standard inner product. Given a set partition $\pi$
of $[2b]$, define the following endomorphism $P_{\pi}$ of $(\C^{n})^{\otimes b}$
by 
\[
\langle P_{\pi}(e_{i_{1}}\otimes\cdots\otimes e_{i_{b}}),e_{i_{b+1}}\otimes\cdots\otimes e_{i_{2b}}\rangle\eqdf\begin{cases}
1 & \text{if \ensuremath{i_{p}=i_{q}} iff \ensuremath{p} and \ensuremath{q} in same block of \ensuremath{\pi}}\\
0 & \text{else.}
\end{cases}
\]
Denote by $\Part([2b])$ the set partitions of $[2b]$, and by $\iota$
the inclusion $S_{b}\hookrightarrow\Part([2b])$ defined by viewing
an element of $S_{b}$ as a matching between two horizontal rows of
$b$ elements, the first row labeled by $[b]$ and the second row
labeled by $[2b]\backslash[b]$. Say $\pi_{1}\leq\pi_{2}$ if the
blocks of $\pi_{1}$ subdivide those of $\pi_{2}$. Let $\subperm(b)\subseteq \Part([2b])$ denote ``sub-permutations'': set partitions $\pi$ of $[2b]$ such that $\pi\leq\iota(\sigma)$ for some $\sigma\in S_b$, namely, such that every block of $\pi$ contains at most one element from each of the two horizontal rows. We write $|\pi|$ for the number of blocks of $\pi$. Here is Cassidy's formula.\footnote{We point out that we do not use the
full force of Cassidy's work in this paper, and it can be bypassed by the (stable) character theory of $S_{n}$. In particular, all we need here is to know that for large $n$, $\chi_{\lambda^{+}(n)}(g)$
is a fixed polynomial of $\fix(g^{k})$ for finitely many $k$ (depending
on $\lambda$) where the polynomial has rational coefficients of $n$
with denominators that are `not too bad'. The more classical formula
\cite[(B.1)]{zbMATH07693381} would accomplish the same thing.
Using Cassidy's result is handy here, however, and we believe that
this formalism will also be useful for some future investigations
in this area.}

\begin{thm}[{Cassidy \cite[Thm.~1.1]{cassidy2023projectionformulasrefinementschurweyljones}}]
\label{thm:Cassidy}For $b\in\Z_+$, $n\geq2b$, and $\lambda\vdash b$,
\begin{equation}
\p_{\lambda}=(-1)^{b}d_{\lambda^{+}(n)}\sum_{\pi\in\subperm(b)}\frac{(-1)^{|\pi|}}{(n)_{|\pi|}}\left(\sum_{\substack{\tau\in S_{b}\\
\pi\leq\iota(\tau)
}
}\chi_{\lambda}(\tau)\right)P_{\pi}\label{eq:projection-formula}
\end{equation}
 is the orthogonal projection in $(\C^{n})^{\otimes b}$ onto $d_{\lambda}$
copies of $V^{\lambda^{+}(n)}$ as an $S_{n}$-representation.
\end{thm}

\begin{proof}[Proof of Proposition \ref{prop:rationality}]
Here we evaluate $\Theta_{\lambda^{+}(n)}\left(Y\right)$ differently
to \cite{MPasympcover}. If 
\[
\rho_{b}:S_{n}\to\End\left(\left(\C^{n}\right)^{\otimes b}\right)
\]
 is the $b$\textsuperscript{th} tensor power of the standard permutation
representation of $S_{n}$, we can rewrite (\ref{eq:theta-first})
as
\[
\Theta_{\lambda^{+}(n)}\left(Y\right)=\frac{\prod_{f\in a,b,c,d}(n)_{\e_{f}}}{d_{\lambda}}\int_{g_{f}\in S_{n}}\mathbf{1}\{g_{f}\text{ obey \,\ensuremath{Y})\}}{\Tr}_{(\C^{n})^{\otimes b}}\left(\p_{\lambda}\rho_{b}\left([g_{a},g_{b}][g_{c},g_{d}]\right)\right)
\]
where $\p_{\lambda}$ is the projection operator from Theorem \ref{thm:Cassidy}. 

From (\ref{eq:projection-formula}), $\p_{\lambda}$ is a linear combination
of partition operators $P_{\pi};$ all the coefficients are rational
functions of $n$ in $n\geq2b$ and can be put over a common denominator
of $(n)_{2b}$. Here we used that $d_{\lambda^{+}(n)}$ is a polynomial
of $n$ in this range by the hook-length formula \cite{Frame_Robinson_Thrall_1954}.
From this we learn $\Theta_{\lambda^{+}(n)}\left(Y\right)$ is a linear
combination of 
\begin{equation}
\int_{g_{f}\in S_{n}}\mathbf{1}\{g_{f}\text{ obey \,\ensuremath{Y}\}}{\Tr}_{(\C^{n})^{\otimes b}}\left(P_{\pi}\rho_{b}\left([g_{a},g_{b}][g_{c},g_{d}]\right)\right)\label{eq:reduce-to-partition}
\end{equation}
with rational coefficients in $n$, with common denominator $(n)_{2b}$.

We can now expand 
\begin{multline}
\label{eq:bigtraceexpansion}
{\Tr}_{(\C^{n})^{\otimes b}}\left(P_{\pi}\rho_{b}\left([g_{a},g_{b}][g_{c},g_{d}]\right)\right)  = \\\sum_{I,J,K,L,M,N,O,Q,R}
 (P_{\pi})_{RI}(g_{a})_{IJ}(g_{b})_{JK}(g_{a})_{LK}(g_{b})_{ML}(g_{c})_{MN}(g_{d})_{NO}(g_{c})_{QO}(g_{d})_{RQ}
\end{multline}
with $I,J,K,L,M,N,O,Q,R\in[n]^{b}$, where by slight abuse of notation
we write $(g_{f})_{IJ}\eqdf\rho_{b}(g_{f})_{IJ}$. We view $I,J,K,L,M,N,O,Q,R$
as well as the fixed labeling of $Y^{(0)}$ concatenated together
as a function 
\[
\Omega:Y^{(0)}\cup[9b]\to[n].
\]
The function $\Omega$ and all these indices contain exactly the same information. We now partition the summation in (\ref{eq:bigtraceexpansion}) according to the partition $p_{\Omega}$ of $Y^{(0)}\cup[9b]$ that the function $\Omega$ defines: two locations are in the same block of $p_{\Omega}$ exactly when they have the same image. Write $\pi\to\Omega$ if and only if the values of $R,I$ encoded by $\Omega$ satisfy $(P_{\pi})_{RI}=1$. We obtain 
\begin{align}
 & \mathbf{1}\{g_{f}\text{ obey \,\ensuremath{Y}\}}{\Tr}_{(\C^{n})^{\otimes b_{}}}\left(P_{\pi}\rho_{b}\left([g_{a},g_{b}][g_{c},g_{d}]\right)\right)\nonumber \\
 & =\sum_{\pi\to\Omega}\mathbf{1}\{g_{f}\text{ obey \,\ensuremath{Y}\}}(g_{a})_{IJ}(g_{b})_{JK}(g_{a})_{LK}(g_{b})_{ML}(g_{c})_{MN}(g_{d})_{NO}(g_{c})_{QO}(g_{d})_{RQ}.\label{eq:summand}
\end{align}
In the above sum the labeling of $Y^{(0)}$ is fixed as usual. Now we use the following elementary integration formula:
\begin{align*}
 & \int_{g\in S_{n}}g_{i_{1}j_{1}}g_{i_{2}j_{2}}\cdots g_{i_{k}j_{k}}\\
 & =\begin{cases}
\frac{1}{(n)_{|\theta|}} & \text{if the maps \ensuremath{\ell\mapsto i_{\ell}} and \ensuremath{\ell\mapsto j_{\ell}} define the same partition \ensuremath{\theta\in\Part([k])}}\\
0 & \text{otherwise.}
\end{cases}
\end{align*}
Notice that $\mathbf{1}\{g_{f}\text{ obey \,\ensuremath{Y}\}}$ is a product of matrix entries of the $g_{f}$:
\[
\mathbf{1}\{g_{f}\text{ obey \,\ensuremath{Y}\}}=\prod_{f\in\{a,b,c,d\}}\prod_{\{(i,j)\,:\,Y\text{ has \ensuremath{f}-colored edge from \ensuremath{i\to j}\}}}(g_{f})_{ij}.
\]
Hence for $n\geq\v+9b$ the integral of the summand of (\ref{eq:summand}) depends only on $p_{\Omega}$ and is either zero or equal to  
\begin{equation}
\frac{1}{\prod_{f\in\{a,b,c,d\}}(n)_{|(p_{\Omega})_{f}|}}\label{eq:integral-summand}
\end{equation}
 where $(p_{\Omega})_{f}$ is the partition induced by $p_{\Omega}$
on 
\begin{itemize}
\item vertices of $Y$ with an outgoing $f$-colored edge, together with
\item elements of $[9b]\backslash Y^{(0)}$ that are in the domain of indices
at the start of the matrix coefficient, e.g.~for $f=a$, domain of
$I$ or domain of $L$. 
\end{itemize}
Note that each 
\[
|(p_{\Omega})_{f}|\leq\v+2b.
\]
Therefore the total contribution from a fixed $p_{\Omega}$ to (\ref{eq:reduce-to-partition})
is (\ref{eq:integral-summand}) times the number of $\Omega$ inducing
$p_{\Omega}$. If a block of $p_{\Omega}$ contains an element of
$Y^{(0)}$ then the corresponding value of $\Omega$ is fixed by the
fixed labeling of $Y^{(0)}$. The other blocks have values in $[n]$
that must be disjoint from each other, hence the number of $\Omega$
for a given $p_{\Omega}$ is 
\[
(n)_{\text{num. blocks of \ensuremath{p_{\Omega}} disjoint from \ensuremath{Y^{(0)}}}}.
\]
Putting the previous arguments together proves the result.

(For general $g$ the only change is the number of factors $(n)_{|(p_{\Omega})_{f}|}$
in the denominator and the end result is a denominator $(n)_{\v+2b}^{1+2g}$.)
\end{proof}

\subsection{Proof of Theorem \ref{thm:asymp-effective}}
\label{sec:proofeffpoly}

We can now complete the proof.

\begin{proof}[Proof of Theorem \ref{thm:asymp-effective}]
Suppose that $|\gamma|\leq q$. We write $\zeta(n)=\zeta(2;S_{n})$.
Now by (\ref{eq:resolution-bound}), (\ref{eq:resolution-expansion}),
and Proposition \ref{prop:tail-cut}, for $n\geq28q^{2}$
\begin{align*}
\E_{g,n}[\fix_{\gamma}] & =\frac{1}{\zeta(n)}\sum_{r\in\mathcal{R}}\frac{(n)_{\v(W_{r})}}{\prod_{f\in a,b,c,d}(n)_{\e_{f}(W_{r})}}\sum_{\substack{\lambda\vdash n\\
\min(b_{\lambda},b_{\lambda}^{\vee})\leq4q
}
}d_{\lambda}\Theta_{\lambda}\left(W_{r}\right)\\
 & \,+O\left((Cq)^{Cq}n^{-q}\right).
\end{align*}
Since above either $b_{\lambda}\leq4q$ or $b_{\lambda}^{\vee}\leq4q$,
for $n\geq Cq$ these events cannot happen simultaneously. Moreover,
it is easy to check from (\ref{eq:theta-first}) that $d_{\lambda}\Theta_{\lambda}\left(W_{r}\right)$
is invariant under $\lambda\mapsto\lambda^{\vee}$ and hence our main
term can be written
\[
\frac{2}{\zeta(n)}\sum_{r\in\mathcal{R}}\frac{(n)_{\v(W_{r})}}{\prod_{f\in a,b,c,d}(n)_{\e_{f}(W_{r})}}\sum_{|\lambda|\leq4q}d_{\lambda^{+}(n)}\Theta_{\lambda^{+}(n)}\left(W_{r}\right).
\]
By the hook-length formula \cite{Frame_Robinson_Thrall_1954} each
$d_{\lambda^{+}(n)}$ is a divisor in $\Q[n]$ of 
\[
(n)_{8q}=n(n-1)\cdots(n-8q+1)
\]
and in particular, a polynomial of $n$ in the range we consider.
Since for everything above, $\v(W_{r}),\e_{f}(W_{r})\leq q,$ for
$n\geq Cq$ Proposition \ref{prop:rationality} tells us that
\[
\Phi_{\gamma}(n)\eqdf\sum_{r\in\mathcal{R}}\frac{(n)_{\v(W_{r})}}{\prod_{f\in a,b,c,d}(n)_{\e_{f}(W_{r})}}\sum_{|\lambda|\leq4q}d_{\lambda^{+}(n)}\Theta_{\lambda^{+}(n)}\left(W_{r}\right)
\]
is a rational function of $n$ whose denominator can be taken to be
$(n)_{9q}^{9}.$ (For general $g$, this is replaced by $(n)_{9q}^{1+4g}$.)

So we write 
\[
\Phi_{\gamma}(n)=\frac{p(n)}{(n)_{9q}^{9}}=\frac{1}{2}\zeta(n)\E_{g,n}[\fix_{\gamma}]+O\left((Cq)^{Cq}n^{-q}\right)
\]
using $\zeta\to2$ as $n\to\infty$ (see (\ref{eq:zeta-limit})) to
remove $\zeta$ from the error. Using the same fact again with a priori
bound $|\E_{g,n}[\fix_{\gamma}]|\leq n$ we learn that $|\frac{1}{n}\frac{p(n)}{(n)_{9q}^{9}.}|$
is bounded as $n\to\infty$ and hence 
\[
\deg(p)\leq81q+1.
\]

Let
\[
g_{q}(t)=\prod_{k=0}^{9q-1}\left(1-kt\right)^{9}
\]
and for $t\eqdf n^{-1}$ write
\[
\Phi_{\gamma}(n)=\frac{p(n)}{(n)_{9q}^{9}}=\frac{n^{\deg(p)}}{n^{81q}}\frac{P(t)}{g_{q}(t)}=t^{81q-\deg(p)}\frac{P(t)}{g_{q}(t)}\eqdf\frac{1}{t}\frac{Q(t)}{g_{q}(t)}.
\]
where $P(t)$ is a polynomial of $t$ of degree $\leq\deg(p)$. For
$\tau\in[0,\frac{1}{Cq^{2}}]$ we have 
\begin{equation}
g_{q}(\tau)\geq C^{-1}\label{eq:g_bound}
\end{equation}
 so we have 
\begin{equation}
Q(t)=t^{81q-\deg(p)+1}P(t)=O(1)\label{eq:Q_upper}
\end{equation}
for $t\in\left[0,\frac{1}{Cq^{2}}\right]\cap\N^{-1}$ (using the a
priori bound again). 
Note that $\deg(Q)\le 81q+1$.
(By the previous remarks, for general $g$, $C$ can be taken to be
linear in $g$.)

To conclude, we begin with
\begin{equation}
\E_{g,n}[\fix_{\gamma}]=\frac{2}{t\zeta(n)}\frac{Q(t)}{g_{q}(t)}+O\left((Cq)^{Cq}n^{-q}\right).\label{eq:temp}
\end{equation}
We now need to deal with the rogue $\zeta(n)$ factor. We can use
a sort of trick to do this efficiently. If $\gamma=\id$ then all
the previous arguments apply to give
\[
t^{-1}=n=\E_{g,n}[\fix_{\id}]=\frac{2}{\zeta(n)}t^{-1}\frac{Q_{\id}(t)}{g_{q}(t)}+O\left((Cq)^{Cq}n^{-q}\right)
\]
where $Q_{\id}$ is a polynomial that also satisfies (\ref{eq:Q_upper}).
Rearranging, 
\begin{align}
\frac{Q_{\id}(t)}{g_{q}(t)} & =\zeta(n)\left(\frac{1}{2}+O\left((Cq)^{Cq}n^{-q-1}\right)\right).\label{eq:Qid_est}
\end{align}
So when $n\geq(Cq)^{2C}$ we have 
\[
\frac{2}{\zeta(n)}=\frac{g_{q}(t)}{Q_{\id}(t)}\left(1+O\left((Cq)^{Cq}n^{-q-1}\right)\right).
\]
Therefore 
\[
\E_{g,n}[\fix_{\gamma}]=t^{-1}\frac{Q(t)}{Q_{\id}(t)}\left(1+O\left((Cq)^{Cq}n^{-q-1}\right)\right)+O\left((Cq)^{Cq}n^{-q}\right).
\]
Note that $Q_{\id}\geq C^{-1}$ in $n\geq q^{C}$ from (\ref{eq:Qid_est})
and (\ref{eq:g_bound}). So by the bound on $Q$ from (\ref{eq:Q_upper})
(also for $Q_{\id}$) we obtain from Lemma \ref{lem:rational-final}
that for $\tau\in[0,\frac{1}{Cq^{2C}}]$
\[
\left|\left(\frac{Q}{Q_{\id}}\right)^{(i)}(\tau)\right|\leq C(Cq)^{Ci}.
\]
Therefore by Taylor's theorem there are $\a_{-1},\a_{0},\ldots,\a_{q}$
(depending on $\gamma)$ with
\[
\left(\frac{Q}{Q_{\id}}\right)(\tau)=\a_{-1}+\a_{0}\tau+\cdots+\a_{q-1}\tau^{q}+O\left((Cq)^{Cq}\tau^{q+1}\right)
\]
for $\tau\in[0,\frac{1}{Cq^{2C}}]$. Thus for $n\geq Cq^{C}$
\[
\E_{g,n}[\fix_{\gamma}]=\a_{-1}n+\a_{0}+\a_{1}n^{-1}+\cdots+\a_{q-1}n^{-(q-1)}+O\left((Cq)^{Cq}n^{-q}\right).
\]
This concludes the proof.
\end{proof}

\section{Geometry of proper powers\label{sec:Proper-powers}}

The aim of this section is to prove Theorem \ref{thm:powers}.

\subsection{Embedding the Cayley graph of $\Gamma$ in $\protect\H$}
\label{sec:embedding}

Denote by $\cay$ the Cayley graph of $\G$ with respect to the 
generators $\{a_{1},b_{1},\ldots,a_{g},b_{g}\}$ in \eqref{eq:presentation of Gamma}. We use the following standard embedding of $\cay$ in the hyperbolic plane $\H$, which is useful, inter alia, as it agrees with the description in \cite{BirmanSeries}. 

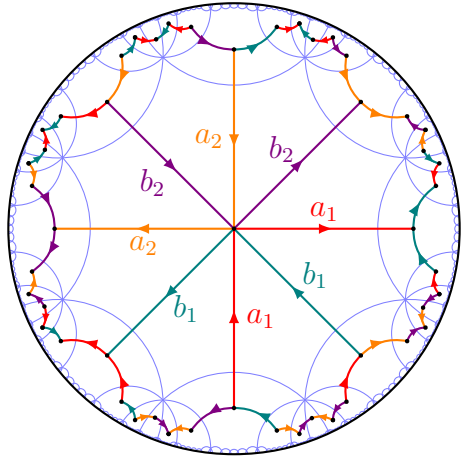
\begin{figure}
\centering \begin{tikzpicture}[scale=3]


\begin{scope}

\clip (0,0) circle (1);

\foreach \i in {1,...,8}
{
	\hgline{45*(\i-1)-25}{45*(\i-1)+25};
	\hgline{45*(\i-1)-30}{45*(\i-1)-6};
	\hgline{45*(\i-1)-15}{45*(\i-1)-39};

	\hgline{45*(\i-1)-3}{45*(\i-1)-6.5};	
	\hgline{45*(\i-1)-.75}{45*(\i-1)-3.5};	
	\hgline{45*(\i-1)+1.25}{45*(\i-1)-1.25};	
	\hgline{45*(\i-1)+3}{45*(\i-1)+6.5};	
	\hgline{45*(\i-1)+.75}{45*(\i-1)+3.5};	

	\hgline{45*(\i-1)-8.5}{45*(\i-1)-5.5};	
	\hgline{45*(\i-1)-9.6}{45*(\i-1)-8.2};	
	\hgline{45*(\i-1)-11}{45*(\i-1)-9.6};	
	\hgline{45*(\i-1)-12.4}{45*(\i-1)-11};	
	\hgline{45*(\i-1)-13.8}{45*(\i-1)-12.4};	
	\hgline{45*(\i-1)-13.5}{45*(\i-1)-15.5};	

	\hgline{45*(\i-1)+8.5}{45*(\i-1)+5.5};	
	\hgline{45*(\i-1)+9.6}{45*(\i-1)+8.2};	
	\hgline{45*(\i-1)+11}{45*(\i-1)+9.6};	
	\hgline{45*(\i-1)+12.4}{45*(\i-1)+11};	
	\hgline{45*(\i-1)+13.8}{45*(\i-1)+12.4};	
	\hgline{45*(\i-1)+13.5}{45*(\i-1)+15.5};	

	\hgline{45*(\i-1)+16.5}{45*(\i-1)+14.5};	
	\hgline{45*(\i-1)-16.5}{45*(\i-1)-14.5};	

	\hgline{45*(\i-1)-26}{45*(\i-1)-24.5};	
	\hgline{45*(\i-1)-25.5}{45*(\i-1)-24};	
	\hgline{45*(\i-1)+26}{45*(\i-1)+24.5};	
	\hgline{45*(\i-1)+25.5}{45*(\i-1)+24};	

	\hgline{45*(\i-1)+23.3}{45*(\i-1)+24.1};	
	\hgline{45*(\i-1)+22.5}{45*(\i-1)+23.3};	
	\hgline{45*(\i-1)+21.7}{45*(\i-1)+22.5};
	\hgline{45*(\i-1)+20.9}{45*(\i-1)+21.7};

	\hgline{45*(\i-1)-17}{45*(\i-1)-16.25};	
	\hgline{45*(\i-1)-17.75}{45*(\i-1)-17};	
	\hgline{45*(\i-1)-18.5}{45*(\i-1)-17.75};	
	\hgline{45*(\i-1)-19.25}{45*(\i-1)-18.5};	
	\hgline{45*(\i-1)+17}{45*(\i-1)+16.25};	
	\hgline{45*(\i-1)+17.75}{45*(\i-1)+17};	
	\hgline{45*(\i-1)+18.5}{45*(\i-1)+17.75};	
	\hgline{45*(\i-1)+19.25}{45*(\i-1)+18.5};	

}

\end{scope}


\draw[thick,red,myptr] (0,0) -- (0.795,0);
\draw[red] (0.4,0.075) node {$a_1$};
\draw[thick,teal,myptr] (0.795,0) to[out=90,in=215] (0.895,0.19);
\draw[thick,teal,myptr] (0.895,-0.19) to[in=-90,out=-215] (0.795,0);
\draw[thick,red,myptr2] (0.91,0.29) to[out=225,in=115] (0.895,0.19);
\draw[thick,red,myptr2b] (0.91,-0.29) to[out=-225,in=-115] (0.895,-0.19);
\draw[thick,teal,myptr2] (0.84,0.345) to[out=300,in=170] (0.91,0.29);
\draw[thick,orange,myptr2b] (0.84,0.345) to[out=110,in=240] (0.85,0.44);
\fill (0.795,0) circle (0.01);
\fill (0.895,0.19) circle (0.01);
\fill (0.895,-0.19) circle (0.01);
\fill (0.91,0.29) circle (0.01);
\fill (0.91,-0.29) circle (0.01);
\fill (0.84,0.345) circle (0.01);

\begin{scope}[rotate=-45]
\draw[thick,teal,myptr] (0.795,0) -- (0,0);
\draw[teal] (0.4,0.11) node {$b_1$};
\draw[thick,orange,myptr] (0.795,0) to[out=90,in=215] (0.895,0.19);
\draw[thick,red,myptr] (0.895,-0.19) to[in=-90,out=-215] (0.795,0);
\draw[thick,violet,myptr2b] (0.895,0.19) to[in=225,out=115] (0.91,0.29);
\draw[thick,violet,myptr2] (0.895,-0.19) to[in=-225,out=-115] (0.91,-0.29);
\draw[thick,orange,myptr2] (0.84,0.345) to[out=300,in=170] (0.91,0.29);
\draw[thick,violet,myptr2] (0.85,0.44) to[in=110,out=240] (0.84,0.345); 
\fill (0.895,0.19) circle (0.01);
\fill (0.895,-0.19) circle (0.01);
\fill (0.795,0) circle (0.01);
\fill (0.91,0.29) circle (0.01);
\fill (0.91,-0.29) circle (0.01);
\fill (0.84,0.345) circle (0.01);
\end{scope}

\begin{scope}[rotate=-90]
\draw[thick,red,myptr] (0.795,0) -- (0,0);
\draw[red] (0.4,0.12) node {$a_1$};
\draw[thick,violet,myptr] (0.795,0) to[out=-90,in=-215] (0.895,-0.19);
\draw[thick,teal,myptr] (0.895,0.19) to[in=90,out=215] (0.795,0);
\draw[thick,orange,myptr2b] (0.895,0.19) to[in=225,out=115] (0.91,0.29);
\draw[thick,orange,myptr2] (0.895,-0.19) to[in=-225,out=-115] (0.91,-0.29);
\draw[thick,violet,myptr2b] (0.91,0.29) to[in=300,out=170] (0.84,0.345);
\draw[thick,orange,myptr2] (0.85,0.44) to[in=110,out=240] (0.84,0.345); 
\fill (0.895,0.19) circle (0.01);
\fill (0.895,-0.19) circle (0.01);
\fill (0.795,0) circle (0.01);
\fill (0.91,0.29) circle (0.01);
\fill (0.91,-0.29) circle (0.01);
\fill (0.84,0.345) circle (0.01);
\end{scope}

\begin{scope}[rotate=-135]
\draw[thick,teal,myptr] (0,0) -- (0.795,0);
\draw[teal] (0.4,0.1) node {$b_1$};
\draw[thick,red,myptr] (0.795,0) to[out=-90,in=-215] (0.895,-0.19);
\draw[thick,red,myptr] (0.895,0.19) to[in=90,out=215] (0.795,0);
\draw[thick,teal,myptr2] (0.91,0.29) to[out=225,in=115] (0.895,0.19);
\draw[thick,teal,myptr2b] (0.91,-0.29) to[out=-225,in=-115] (0.895,-0.19);
\draw[thick,orange,myptr2b] (0.91,0.29) to[in=300,out=170] (0.84,0.345);
\draw[thick,violet,myptr2b] (0.84,0.345) to[out=110,in=240] (0.85,0.44);
\fill (0.895,0.19) circle (0.01);
\fill (0.895,-0.19) circle (0.01);
\fill (0.795,0) circle (0.01);
\fill (0.91,0.29) circle (0.01);
\fill (0.91,-0.29) circle (0.01);
\fill (0.84,0.345) circle (0.01);
\end{scope}

\begin{scope}[rotate=-180]
\draw[thick,orange,myptr] (0,0) -- (0.795,0);
\draw[orange] (0.4,0.075) node {$a_2$};
\draw[thick,violet,myptr] (0.795,0) to[out=90,in=215] (0.895,0.19);
\draw[thick,violet,myptr] (0.895,-0.19) to[in=-90,out=-215] (0.795,0);
\draw[thick,orange,myptr2] (0.91,0.29) to[out=225,in=115] (0.895,0.19);
\draw[thick,orange,myptr2b] (0.91,-0.29) to[out=-225,in=-115] (0.895,-0.19);
\draw[thick,violet,myptr2] (0.84,0.345) to[out=300,in=170] (0.91,0.29);
\draw[thick,red,myptr2b] (0.84,0.345) to[out=110,in=240] (0.85,0.44);
\fill (0.795,0) circle (0.01);
\fill (0.895,0.19) circle (0.01);
\fill (0.895,-0.19) circle (0.01);
\fill (0.91,0.29) circle (0.01);
\fill (0.91,-0.29) circle (0.01);
\fill (0.84,0.345) circle (0.01);
\end{scope}

\begin{scope}[rotate=-225]
\draw[thick,violet,myptr] (0.795,0) -- (0,0);
\draw[violet] (0.4,0.11) node {$b_2$};
\draw[thick,red,myptr] (0.795,0) to[out=90,in=215] (0.895,0.19);
\draw[thick,orange,myptr] (0.895,-0.19) to[in=-90,out=-215] (0.795,0);
\draw[thick,teal,myptr2b] (0.895,0.19) to[in=225,out=115] (0.91,0.29);
\draw[thick,teal,myptr2] (0.895,-0.19) to[in=-225,out=-115] (0.91,-0.29);
\draw[thick,red,myptr2] (0.84,0.345) to[out=300,in=170] (0.91,0.29);
\draw[thick,teal,myptr2] (0.85,0.44) to[in=110,out=240] (0.84,0.345); 
\fill (0.895,0.19) circle (0.01);
\fill (0.895,-0.19) circle (0.01);
\fill (0.795,0) circle (0.01);
\fill (0.91,0.29) circle (0.01);
\fill (0.91,-0.29) circle (0.01);
\fill (0.84,0.345) circle (0.01);
\end{scope}

\begin{scope}[rotate=-270]
\draw[thick,orange,myptr] (0.795,0) -- (0,0);
\draw[orange] (0.4,0.11) node {$a_2$};
\draw[thick,teal,myptr] (0.795,0) to[out=-90,in=-215] (0.895,-0.19);
\draw[thick,violet,myptr] (0.895,0.19) to[in=90,out=215] (0.795,0);
\draw[thick,red,myptr2b] (0.895,0.19) to[in=225,out=115] (0.91,0.29);
\draw[thick,red,myptr2] (0.895,-0.19) to[in=-225,out=-115] (0.91,-0.29);
\draw[thick,teal,myptr2b] (0.91,0.29) to[in=300,out=170] (0.84,0.345);
\draw[thick,red,myptr2] (0.85,0.44) to[in=110,out=240] (0.84,0.345); 
\fill (0.895,0.19) circle (0.01);
\fill (0.895,-0.19) circle (0.01);
\fill (0.795,0) circle (0.01);
\fill (0.91,0.29) circle (0.01);
\fill (0.91,-0.29) circle (0.01);
\fill (0.84,0.345) circle (0.01);
\end{scope}

\begin{scope}[rotate=45]
\draw[thick,violet,myptr] (0,0) -- (0.795,0);
\draw[violet] (0.4,0.1) node {$b_2$};
\draw[thick,orange,myptr] (0.795,0) to[out=-90,in=-215] (0.895,-0.19);
\draw[thick,orange,myptr] (0.895,0.19) to[in=90,out=215] (0.795,0);
\draw[thick,violet,myptr2] (0.91,0.29) to[out=225,in=115] (0.895,0.19);
\draw[thick,violet,myptr2b] (0.91,-0.29) to[out=-225,in=-115] (0.895,-0.19);
\draw[thick,red,myptr2b] (0.91,0.29) to[in=300,out=170] (0.84,0.345);
\draw[thick,teal,myptr2b] (0.84,0.345) to[out=110,in=240] (0.85,0.44);
\fill (0.895,0.19) circle (0.01);
\fill (0.895,-0.19) circle (0.01);
\fill (0.795,0) circle (0.01);
\fill (0.91,0.29) circle (0.01);
\fill (0.91,-0.29) circle (0.01);
\fill (0.84,0.345) circle (0.01);
\end{scope}

\foreach \i in {0,45,90,135,180,225,270,-45}
{
\begin{scope}[rotate=\i]
\fill (0.91,-0.29) circle (0.01);
\end{scope}
}

\fill (0,0) circle (0.01);

\draw[thick] (0,0) circle (1);

\end{tikzpicture} \caption{\label{fig:tiling}Illustration of the tiling $\mathcal{T}$ and embedded Cayley graph of $\Gamma$ for $g=2$.}
\end{figure}

Consider the tiling ${\cal T}$ of $\H$ by regular $4g$-gons with geodesic sides and interior angles $\frac{2\pi}{4g}$. The Cayley graph is then the tiling dual to ${\cal T}$. More precisely, the vertices of $\cay$ are the centers of the $4g$-gons, and the edges are geodesic arcs between the centers of any two $4g$-gons sharing a side. Every edge of $\cay$ is directed and labeled by one of the $2g$ generators. For any generator $x$, when traversing an $x$-edge against its direction, one reads $x^{-1}$. At each vertex, there is exactly one edge directed outward and one directed inward with every given label. The cyclic order of the $4g$ outgoing edges, say clockwise, is
\[
a_{1},b_{1}^{-1},a_{1}^{-1},b_{1},a_{2},b_{2}^{-1},a_{2}^{-1},b_{2},\ldots,b_{g}.
\]
(Up to isometry, there is basically a unique way to label the edges
of the dual tiling in this way in a compatible fashion.) The boundary
of every dual $4g$-gon now reads cyclically, counterclockwise, the
relation $[a_{1},b_{1}]\cdots[a_{g},b_{g}]$. See Figure \ref{fig:tiling}.

We pick one of the vertices to be the identity element and mark it
by $o$. Any other vertex $v$ of $\cay$ now corresponds to the element
of $\G$ represented by any of the paths from $o$ to $v$. The left
action of $\G$ on $\cay$ extends to an action by isometries on $\H$.
The isometry corresponding to each non-trivial element of $\G$ is
hyperbolic, which means that it admits a (unique,
bi-infinite) geodesic axis in $\H$: the isometry acts on this axis
by translation. These facts can be found, e.g., in \cite[\S1]{farb2011primer}.

\subsection{$\prod$-shaped paths for elements in $\protect\G$ }

The following definition will play a key role in the sequel. Roughly speaking, it is the appropriate replacement for surface groups of the fact that any element $\gamma$ of a free group
can be expressed as $bhb^{-1}$ where $h$ is the cyclic reduction of $\gamma$.

\begin{defn}[$\prod$-path]
\label{def:prod-path} Let $1\ne\gamma\in\G$. Consider the unique
(bi-infinite) geodesic $\rho\subset\H$ which is the axis of the isometry
given by $\gamma$. The $\prod$-path corresponding to $\gamma$,
which we denote by $\prod_{\gamma}$, is the path in $\H$ made out
of three geodesics arcs: a geodesic arc $[o,x]$ from $o$ to $x\in\rho$
which is perpendicular to $\rho$, the geodesic arc $[x,\gamma.x]\subset\rho$, and the geodesic arc $[\gamma.x,\gamma.o]$ (which is, too, perpendicular to $\rho$).
\end{defn}

See Figure \ref{fig:chet} below. Note that for any hyperbolic isometry
of $\H$, each half-space at either side of the axis is invariant
under the isometry. This means that $o$ and $\gamma.o$ lie on the
same side of the axis, so that $\prod_{\gamma}$ is made of three
geodesic arcs with two right-turns, or two left-turns, between them
(unless $o$ happens to lie on the axis of $\gamma$, in which case
$\prod_{\gamma}$ is simply a geodesic arc). 

We will need the following easy fact from hyperbolic geometry:
\begin{fact}
\label{fact:about Chet-shape}Let $L$ be a geodesic in $\H$ and
$L_{1},L_{2}$ be two geodesic rays perpendicular to $L$, emanating
from two distinct points $y_{1}$ and $y_{2}$, respectively, on $L$,
and at the same side of $L$. Denote $d_{1}=d(y_{1},y_{2})>0$. Let
$z$ be a point on $L_{2}$ at distance $d_{2}\ge0$ from $y_{2}$.
Then the distance $d(z,L_{1})$ of $z$ from $L_{1}$ grows monotonically
with each of $d_{1}$ and $d_{2}$, and for every fixed $d_{1}>0$,
\[
d(z,L_{1})\stackrel{d_{2}\to\infty}{\to}\infty.
\]
\end{fact}

\begin{figure}
\centering \begin{tikzpicture}

\draw (0,2.35) -- (0.15,2.35) -- (0.15,2.48);

\begin{scope}[rotate=-40]
\draw (0,2.35) -- (0.15,2.35) -- (0.15,2.48);
\end{scope}

\begin{scope}[yshift=-1.25cm]
\draw (0,2.35) -- (0.15,2.35) -- (0.15,2.48);
\end{scope}

\centerarc[thick,black!10](0,0)(180:0:1.25);

\centerarc[thick,black!35](0,0)(180:0:2.5);

\draw[thick,black!35] (0,0) -- (0,3);

\centerarc[thick,black!35](3.915,0)(180:113:3);

\centerarc[thick,blue](0,0)(90:50:2.5);

\centerarc[thick,blue](3.915,0)(165.5:140.05:3);

\centerarc[thick,dashed,blue!50](0,0)(90:35:1.25);

\draw[thick] (-3,0) -- (3,0);

\draw (2.55,1) node {$L$};

\draw (0.05,3.2) node {$L_1$};

\draw (3,2.8) node {$L_2$};

\filldraw (0,2.5) circle (0.05);
\draw (-.15,2.3) node {$y_1$};

\filldraw (1.61,1.915) circle (0.05);
\draw (1.92,1.93) node {$y_2$};

\filldraw (1.01,0.74) circle (0.05);
\draw (1.2,0.8) node {$z$};

\draw[blue] (0.95,2.55) node {$d_1$};

\draw[blue] (1.45,1.25) node {$d_2$};

\begin{scope}[rotate=-23]
\draw[blue] (0,1.42) node {\rotatebox{-18}{$\scriptstyle d(z,L_1)$}};
\end{scope}

\end{tikzpicture} \caption{\label{fig:Fact 3}Illustration of Fact \ref{fact:about Chet-shape}
(half-plane model of $\mathbb{H}^{2}$)}
\end{figure}
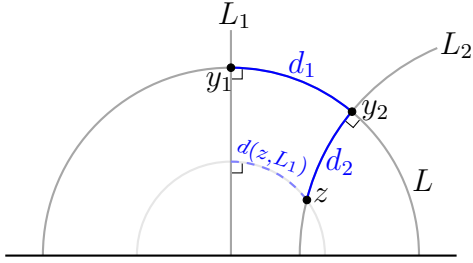

Fact \ref{fact:about Chet-shape} is illustrated in Figure \ref{fig:Fact 3}.
It follows immediately from the explicit formula $\sinh d(z,L_1)=\sinh(d_{1})\cosh(d_{2})$,
that may be found in \cite[Thm.~32.21(a')]{martin1975foundations}.
For completeness, we give an elementary argument.
\begin{proof}[Proof of Fact \ref{fact:about Chet-shape}]
 Assume without loss of generality that $L_{1}$ is the geodesic $x=0$ in the upper half plane model of $\H$. Then the distance between $z$ and $L_{1}$ is determined by the argument of $z$, and monotonically increases and tends to $\infty$ as the argument in $\left[0,\frac{\pi}{2}\right]$
decreases. This proves the claims when $d_{1}>0$ is fixed.

Now fix $d_{2}\ge0$. Then $d(z,L_{1})$ is continuous as a function of $d_{1}$. Clearly $d(z,L_{1})\stackrel{d_{1}\searrow0}{\to}0$
and $d(z,L_{1})\stackrel{d_{1}\to\infty}{\to}\infty$. So if it is not monotonically increasing, there must be two different values $\alpha_{1}<\alpha_{2}$ of $d_{1}$ with two corresponding points $z_{1}\ne z_{2}$ with $d(z_{1},L)=d(z_{2},L)$. Then $z_{1}$ must lie on the Euclidean interval from $z_{2}$ to $0$. So the hyperbolic isometry $\Phi\colon t\mapsto ct$ maps $z_{1}$ to $z_{2}$ for some $c>1$. But this isometry shows that $z_{2}$ is at distance $d_{2}$ from the geodesic $\Phi(L)$ which lies strictly
above $L$ --- a contradiction.
\end{proof}

The following lemma is illustrated in Figure \ref{fig:chet}.

\begin{lem}
\label{lem:in H^2 prod and arc are close}There is a constant $c_{1}=c_{1}(g)>0$ so that $\prod_{\gamma}$ and the geodesic arc $[o,\gamma.o]$ are each contained in a $c_{1}$-neighborhood of the other. Moreover, there are two points $z_{1},z_{2}\in[o,\gamma.o]$, with $z_{2}$ not closer to $o$ than $z_{1}$, such that the two corners $x$ and $\gamma.x$ of $\prod_{\gamma}$ satisfy $d_{\H}(x,z_{1})<c_{1}$
and $d_{\H}(\gamma.x,z_{2})<c_{1}$. 
\end{lem}

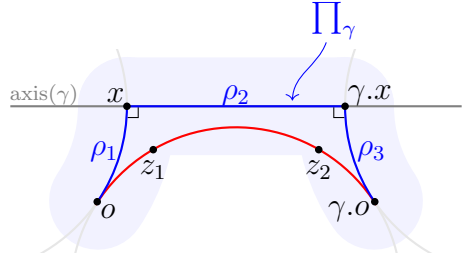
\begin{figure}
\centering \begin{tikzpicture}

\begin{scope}

\clip (-3,0.75) rectangle (3,-1.95);




\centerarc[thick,blue!5,line cap=round,line 
width=1.3cm](-3.7,0)(0.15:-34:2.25);
\centerarc[thick,blue!5,line cap=round,line 
width=1.3cm](3.7,0)(180-0.15:180+34:2.25);
\draw[thick,blue!5,line cap=round,line width=1.3cm] (-1.45,0) -- (1.45,0);


\draw (-1.45,-.15) -- (-1.3,-.15) -- (-1.3,0);
\draw (1.45,-.15) -- (1.3,-.15) -- (1.3,0);

\draw[thick,black!50] (-3,0) -- (3,0);
\draw[thick,black!10] (3.7,0) circle (2.25);
\draw[thick,black!10] (-3.7,0) circle (2.25);


\draw[thick,black!10] (0,-2.48) circle (2.2);
\centerarc[thick,red](0,-2.48)(32:148:2.2);

\centerarc[thick,blue](-3.7,0)(0.15:-34:2.25);
\centerarc[thick,blue](3.7,0)(180-0.15:180+34:2.25);
\draw[thick,blue] (-1.45,0) -- (1.45,0);

\fill (-1.45,0) circle (0.05);
\fill (1.45,0) circle (0.05);
\fill (1.845,-1.27) circle (0.05);
\fill (-1.845,-1.27) circle (0.05);

\draw (-1.7,-1.4) node {$o$};
\draw (1.5,-1.4) node {$\gamma.o$};

\draw (-1.6,0.15) node {$x$};
\draw (1.78,0.15) node {$\gamma.x$};

\draw[black!50] (-2.55,0.15) node {$\scriptstyle \mathrm{axis}(\gamma)$};

\draw[blue] (-1.75,-0.6) node {$\rho_1$};
\draw[blue] (1.8,-0.6) node {$\rho_3$};
\draw[blue] (0,0.15) node {$\rho_2$};

\end{scope}

\draw[blue,->] (1.2,0.85) to[out=250,in=90] (0.75,0.075);
\draw[blue] (1.3,0.75) node[above] {$\prod_\gamma$};


\fill (-1.1,-0.5775) circle (0.05) node[below] {$z_1$};
\fill (1.1,-0.5775) circle (0.05) node[below] {$z_2$};

\end{tikzpicture} \caption{\label{fig:chet}Illustration of Definition \ref{def:prod-path} and
Lemma \ref{lem:in H^2 prod and arc are close}}
\end{figure}

\begin{proof}
$\H$ has $\delta$-thin triangles in the sense
of Rips and Gromov: every side of a geodesic triangle is contained
in a $\delta$-neighborhood of the union of the other two sides \cite[\S III.H.1]{bridson2013metric}.
It follows that every side of a geodesic quadrilateral is contained
in a $2\delta$-neighborhood of the union of the other three sides.
Hence, $[o,\gamma.o]\subset{\cal N}_{2\delta}\left(\prod_{\gamma}\right)$. 

Conversely, denote by $\rho_{1}=[o,x]$, $\rho_{2}=[x,\gamma.x]$
and $\rho_{3}=[\gamma.x,\gamma.o]$ the three geodesic arcs composing
$\prod_{\gamma}$, so $x=\rho_{1}\cap\rho_{2}$ and $\gamma.x=\rho_{2}\cap\rho_{3}$.
Note that $|\rho_{2}|$ is bounded from below by the injectivity radius
of the genus-$g$ hyperbolic surface $\G\backslash\H$. Hence, by
Fact~\ref{fact:about Chet-shape}, there is some $c'\ge2\delta$ such
that if $y\in\rho_{1}$ satisfies $d(y,x)>c'$, then $d(y,\rho_{3})>4\delta$.
If $\left|\rho_{1}\right|<2c'$ we may simply take $z_{1}=o$ and
$z_{2}=\gamma.o$ (recall that $\left|\rho_{3}\right|=\left|\rho_{1}\right|$).
Otherwise, we reason as follows. By construction, the sets ${\cal N}_{2\delta}\left(\left\{ y\in\rho_{1}:d(y,x)>c'\right\} \right)$
and ${\cal N}_{2\delta}(\rho_{3})$ are disjoint. Thus $[o,\gamma.o]$
must visit a point $z_{1}\in{\cal N}_{2\delta}\left(\left\{ y\in\rho_{1}:d(y,x)\le2c'\right\} \right)$
before it reaches ${\cal N}_{2\delta}\left(\rho_{3}\right)$. Clearly,
$d(z_{1},x)\le2c'+2\delta$. Symmetrically, the inverse geodesic $[\gamma.o,o]$
must visit a point $z_{2}\in{\cal N}_{2\delta}\left(\left\{ y\in\rho_{3}:d(y,\gamma.x)\le2c'\right\} \right)$
before it reaches ${\cal N}_{2\delta}\left(\rho_{1}\right)$ and $d(z_{2},\gamma.x)\le2c'+2\delta$.
Finally, 
\[
\rho_{1}\subset{\cal N}_{\delta}\left([o,z_{1}]\cup[z_{1},x]\right)\subseteq{\cal N}_{3\delta+2c'}\left([o,z_{1}]\right),
\]
likewise $\rho_{3}\subset{\cal N}_{3\delta+2c'}\left([z_{2},\gamma.o]\right)$,
and 
\[
\rho_{2}\subset{\cal N}_{2\delta}\left([x,z_{1}]\cup[z_{1},z_{2}]\cup[z_{2},\gamma.x]\right)\subseteq{\cal N}_{4\delta+2c'}\left([z_{1},z_{2}]\right).
\]
The lemma is now proven with $c_{1}=4\delta+2c'$.
\end{proof}

We now follow the terminology from \cite{BirmanSeries}. Recall the
tiling ${\cal T}$ of $\H$ that was introduced in \S\ref{sec:embedding}. Denote by $hO$ the $4g$-gon with center
$h.o$. 

\begin{defn}[{Geodesic edge path \cite[p.~455]{BirmanSeries}}]
\label{def:geod path and word} Let $\rho$ be a finite-length, oriented
geodesic arc in $\H$. Assume it begins in the interior of the $4g$-gon
$h_{1}O$, and then visits $h_{2}O,h_{3}O$ and so on until it terminates
at the interior of $h_{k}O$, and that it does not visit any vertex
of the tiling ${\cal T}$ and does not coincide with any edge of ${\cal T}$.
Then the corresponding \emph{geodesic edge path} is the path $h_{1}.o-h_{2}.o-\ldots-h_{k}.o$
in $\cay$. The corresponding geodesic word is the (reduced) word
in $F(a_{1},\ldots,b_{g})$ 
\[
\left(h_{1}^{-1}h_{2}\right)\left(h_{2}^{-1}h_{3}\right)\ldots\left(h_{k-1}^{-1}h_{k}\right)
\]
(recall that if $hO$ and $h'O$ share a common side, then the edge
in $\cay$ from $hO$ to $h'O$ is labeled by $h^{-1}h'\in\left\{ a_{1}^{\pm1},\ldots,b_{g}^{\pm1}\right\} $).
If $\rho$ goes through a vertex of ${\cal T}$ or coincides with
an edge of ${\cal T}$, we deform it slightly as in {[}Ibid., Fig.~1{]},
and then define the geodesic edge path as above using the deformed
arc.\footnote{\label{fn:deformations}To avoid a vertex, the deformation consists
of taking a detour on either side of the vertex. To avoid an edge,
it consists of moving the entire geodesic arc to either side of the
edge. The geodesic arcs we use here never coincide with edges of ${\cal T}$. If $\rho$ needs to be deformed, the resulting geodesic edge path
is not unique, but this is immaterial for our purposes.}
\end{defn}

Using geodesic edge paths, we have an analog of Lemma \ref{lem:in H^2 prod and arc are close}
in $\cay$. It basically says that in $\cay$ there are, too, a geodesic
path from $o$ to $\gamma.o$ and a $\prod_{\gamma}$-like path which
are close to each other. We focus on proper powers and on the precise
properties that will be needed below. We view $\cay$ as a metric
space with edge-length $1$, and a geodesic path in it is any shortest
path between two vertices (here and in the following lemma, vertices
are vertices of $\cay$ and not of the tiling, unless stated otherwise).
\begin{lem}
\label{lem:breaking proper powers }There is a constant $c_{2}=c_{2}(g)>0$
so that the following holds. Let $k\in\Z_{\ge2}$ and $1\ne\gamma\in\G$
be a proper $k^{\text{th}}$-power. Then there are $b,h\in\G$ with
$\gamma=bh^{k}b^{-1}$, a geodesic path $\tau_{\gamma}\subset\cay$
from $o$ to $\gamma.o$, and four vertices $v_{0},v_{1},v_{2},v_{k}\in\tau_{\gamma}$,
with $v_{j}$ not closer to $o$ than $v_{i}$ when $i<j$, such that
$d_{\cay}(bh^{j}.o,v_{j})<c_{2}$ for all $j=0,1,2,k$.
\end{lem}

\begin{proof}
Let $\gamma_{0}\in\G$ with $\gamma=\gamma_{0}^{\,k}$. Consider the
geodesic arc $\left[o,\gamma.o\right]$ in $\H$, (deform it a bit
if need be), and let $\tau_{\gamma}$ be the associated geodesic edge
path as in Definition \ref{def:geod path and word}. By \cite[Thm.~2.8(b)]{BirmanSeries},
every geodesic edge path, and $\tau_{\gamma}$ in particular, is a
geodesic in $\cay$.

Now consider the $\prod$-path $\prod_{\gamma}$ with vertices $o,x,\gamma.x,\gamma.o$.
Let $b\in\G$ satisfy that $x\in bO$ (if $x$ lies on a vertex or
an edge of ${\cal T}$, pick some neighboring $4g$-gon arbitrarily),
and let $h\eqdf b^{-1}\gamma_{0}b\in\G$. Note that $\gamma=bh^{k}b^{-1}$.
Note also that $\gamma_{0}.x\in\gamma_{0}bO=bhO$, and likewise $\gamma_{0}^{~j}.x\in bh^{j}O$
for all $j\in\Z$. As the axes of $\gamma$ and of $\gamma_{0}$
are identical, the points $x,\gamma_{0}.x,\ldots,\gamma_{0}^{~k}.x=\gamma.x$
all lie on $\left[x,\gamma.x\right]$. By Lemma \ref{lem:in H^2 prod and arc are close},
there are points $z_{0}$ and then $z_{k}$ along $[o,\gamma.o]$
at distance $<c_{1}(g)$ from $x$ and $\gamma.x$, respectively.
Now let $\delta=\delta_{\H}$ be such that triangles in $\H$ are
$\delta$-thin, and consider the geodesic quadrilateral in $\H$ with
corners $x,\gamma.x,z_{k},z_{0}$. The point $\gamma_{0}.x\in[x,\gamma.x]$
is at distance $<2\delta$ from the union of the other three sides,
and thus at distance $<2\delta+c_{1}(g)$ from $[z_{0},z_{k}]$. We
mark a point $z_{1}\in[z_{0},z_{k}]$ with $d_{\H}(z_{1},\gamma_{0}.x)<2\delta+c_{1}(g)$.
Similarly, we find a point $z_{2}\in[z_{1},z_{k}]$ with $d_{\H}(z_{2},\gamma_{0}^{\,2}.x)<4\delta+c_{1}(g)$
(if $k=2$ we may simply set $z_{2}=z_{k}$).

Denote by $D>0$ the diameter of the $4g$-gons. Every point on a
geodesic arc in $\H$ lies in the $D$-neighborhood of its corresponding
geodesic edge path, and even in the $D$-neighborhood of the vertices
of the geodesic edge path. So we may find vertices $v_{0},v_{1},v_{2},v_{k}\in\tau_{\gamma}$
with $d_{\H}(z_{j},v_{j})<D$ for $j=0,1,2,k$. It is clear, by the
definition of a geodesic edge path, that we may choose the $v_{j}$'s
so that $v_{j}$ is not closer than $v_{i}$ to $o$ when $i<j$.
As $d_{\H}(\gamma_{0}^{~j}.x,bh^{j}.o)<D$ for all $j=0,1,2,k$, we
obtain that 
\[
d_{\H}\left(bh^{j}.o,v_{j}\right)<4\delta+c_{1}\left(g\right)+2D.
\]
Finally, by the \u{S}varc-Milnor Lemma (e.g., \cite[Prop.~I.8.19]{bridson2013metric}),
the embedding above of $\cay$ in $\H$ is a quasi-isometry, and bounded
distances remain bounded. In particular, there is some $c_{2}\left(g\right)>0$
so that $d_{\cay}\left(bh^{j}.o,v_{j}\right)<c_{2}\left(g\right)$
for $j=0,1,2,k$.
\end{proof}

The following lemma considers words with letters in an arbitrary (finite)
subset ${\cal S}\subset\G$. Let $w$ be such a word of length $p$.
We let $w_{[t_{1},t_{2})}$ denote the subword consisting of the letters
at positions $t_{1},t_{1}+1,\ldots,t_{2}-1$ (here $1\le t_{1}\le t_{2}\le p+1$),
and $w_{[t_{1},t_{2}]}$ the subword consisting of the letters at
positions $t_{1},\ldots,t_{2}$ (here $1\le t_{1}\le t_{2}\le p$).

\begin{lem}
\label{lem:padding power words to get nice structure}There is a constant
$c_{3}=c_{3}(g)>0$ so that the following holds. Fix any finite subset
$\mathcal{S}\subset\Gamma$ so that $\max_{\gamma\in\mathcal{S}}|\gamma|\le \ell$,
and let $w$ be any not-necessarily-reduced word of length $p$ with
letters in $\mathcal{S}$ representing a non-trivial $k^{\text{th}}$
power in $\Gamma$ for some $k\ge2$. Then there are $1\le t_{1}\le t_{2}\le t_{3}\le t_{4}\le p+1$,
and elements $b,h,u_{1},\ldots,u_{4}\in\Gamma$ with $\left|u_{i}\right|\le c_{3}\left(\ell+\log p\right)$,
such that 
\[
\left(w_{[1,t_{1})}u_{1}\right)\cdot\left(u_{1}^{-1}w_{[t_{1},t_{2})}u_{2}\right)\cdot\left(u_{2}^{-1}w_{[t_{2},t_{3})}u_{3}\right)\cdot\left(u_{3}^{-1}w_{[t_{3},t_{4})}u_{4}\right)\cdot\left(u_{4}^{-1}w_{\left[t_{4},p\right]}\right)
\]
is of the form $b\cdot h\cdot h\cdot h^{k-2}\cdot b^{-1}$. Namely,
\begin{eqnarray*}
 & w_{[1,t_{1})}u_{1}\equiv_{\Gamma}b,\qquad & u_{1}^{-1}w_{[t_{1},t_{2})}u_{2}\equiv_{\Gamma}u_{2}^{-1}w_{[t_{2},t_{3})}u_{3}\equiv_{\Gamma}h\\
 & u_{3}^{-1}w_{[t_{3},t_{4})}u_{4}\equiv_{\Gamma}h^{k-2},\qquad & u_{4}^{-1}w_{\left[t_{4},p\right]}\equiv_{\Gamma}b^{-1}.
\end{eqnarray*}
Here $\equiv_{\Gamma}$ denotes equality in $\Gamma$.
\end{lem}

\begin{proof}
All the distances in this proof are measured in $\cay$. Let $1\ne\gamma\in\G$
denote the element represented by $w$, and let $b,h\in\G$, $\tau_{\gamma}$
and $v_{0},v_{1},v_{2},v_{k}\in\tau_{\gamma}$ be the elements of
$\Gamma$, the geodesic and the vertices given by Lemma \ref{lem:breaking proper powers }.
We may expand $w$ to a word of length at most $\ell p$ in the generators
$\left\{ a_{1}^{\pm},\ldots,b_{g}^{\pm}\right\} $ of $\Gamma$, and
the expanded word then represents a path $\mathfrak{p}_{w}$ of length
at most $\ell p$ from $o$ to $\gamma.o$. The Cayley graph $\cay$ is
a hyperbolic space (in the sense of Rips-Gromov) with $\delta$-thin
triangles for some $\delta>0$. A standard fact in hyperbolic geometry
(e.g., \cite[Prop.~III.H.1.6]{bridson2013metric}) gives that as $\tau_{\gamma}$
is a geodesic in $\cay$ between $o$ and $\gamma.o$, every point
$x\in\tau_{\gamma}$ is at distance at most $\delta\left|\log_{2}\left|\mathfrak{p}_{w}\right|\right|+1\le\delta\left|\log_{2}(\ell p)\right|+1$
from $\mathfrak{p}_{w}$. Together with Lemma \ref{lem:breaking proper powers },
this shows that each of the points $b.o,bh.o,bh^{2}.o$ and $bh^{k}.o$
are rather close to $\mathfrak{p}_{w}$. It remains to show we can
find four points along $\mathfrak{p}_{w}$ which are close to $b.o,bh.o,bh^{2}.o,bh^{k}.o$
and ordered correctly.

First, as every vertex on $\mathfrak{p}_{w}$ is at distance at most
$\ell$ from the beginning of a letter of $w$, there is a vertex $y_{0}$
on $\mathfrak{p}_{w}$ which is at the beginning of a letter of $w$
and with $d(v_{0},y_{0})<\delta\left|\log_{2}\left(\ell p\right)\right|+1+\ell$,
so 
\[
d(b.o,y_{0})\le d(b.o,v_{0})+d(v_{0},y_{0})<c_{2}+\delta\left|\log_{2}\left(\ell p\right)\right|+1+\ell<c'\left(\ell+\log p\right),
\]
for some $c'>0$ depending only on $g$. Set $u_{1}$ to be a word
representing a shortest path in $\cay$ from $y_{0}$ to $b.o$, and
let $t_{1}$ be the smallest integer so that $w_{[1,t_{1})}$ ends
at $y_{0}$. Now consider the path from $v_{0}$ to $\gamma.o$ given
by $\mathfrak{q}\cdot w_{[t_{1},p]}$, where $\mathfrak{q}$ is a
geodesic path from $v_{0}$ to $y_{0}$ (so $\left|\mathfrak{q}\right|<\delta\left|\log_{2}(\ell p)\right|+1+\ell$).
As $v_{1}$ lies on a geodesic from $v_{0}$ to $\gamma.o$, there
is a vertex $y_{1}$ on $\mathfrak{p}_{w}$, at the beginning of a
letter of $w_{[t_{1},p]}$, such that
\[
d(v_{1},y_{1})<\delta\left|\log_{2}\left(\ell(p+1-t_{1})+\left|\mathfrak{q}\right|\right)\right|+1+\ell+\left|\mathfrak{q}\right|,
\]
so $d(bh.o,y_{1})<c''(\ell+\log p)$ for some $c''>0$ depending only
on $g$. We let $u_{2}$ be a word representing a shortest path from
$y_{1}$ to $bh.o$ and $t_{2}$ the smallest integer so that $t_{2}\ge t_{1}$
and $w_{[1,t_{2})}$ ends at $y_{1}$.

In the same manner, and enlarging the constant $c''$ as needed at
each step, we find points $y_{2}$ and $y_{k}$ on the path of $w$,
words $u_{3}$ and $u_{4}$ and times $t_{3}$ and $t_{4}$ so that
the statement of the lemma holds.
\end{proof}

\subsection{Proof of Theorem \ref{thm:powers}}

As is explained in \cite[\S 6]{MageeSalle2}, the proof of Theorem \ref{thm:powers} --- and in fact of Theorem~\ref{thm:completepowers} below, of which Theorem~\ref{thm:powers} is a special case ---
can be reduced to bounding the exponential growth rate of the probability that a product $\boldsymbol{\gamma}_{1}\cdots\boldsymbol{\gamma}_{p}$ of i.i.d.\ finitely supported random group elements $\boldsymbol{\gamma}_i\in\Gamma$ (i.e., a random walk on $\Gamma$ with a general finitely supported generating measure) lands on a proper power.

In the case of the free group, the requisite bound follows by a simple spectral argument that dates back to \cite[Lem.\ 2.4]{friedman2003relative}. This argument can be adapted to the setting of surface groups due to Lemma \ref{lem:padding power words to get nice structure}.

\begin{thm}
\label{thm:completepowers}
For every $d\in\N$ and self-adjoint $x\in\mathrm{M}_d(\C)\otimes\C[\Gamma]$, we have
\[
        \limsup_{p\to\infty} |[\mathrm{tr}\otimes u_1](x^p)|^{1/p} \le
        \|[\mathrm{id}\otimes\lambda](x)\|.
\]
\end{thm}

\begin{proof}
By \cite[Prop.~6.2]{MageeSalle2} and as surface groups
satisfy the rapid decay property \cite[Thm.~3.2.1]{jolissaint1990rapidly}, 
it suffices to prove the theorem in the special case that $x\in \C[\Gamma]$ with positive coefficients.
Let us therefore fix a self-adjoint $x=\sum_{\gamma\in\Gamma}\alpha_{\gamma}\gamma$
with $|x|=\ell$ and $\alpha_\gamma\ge 0$.
By homogeneity, we can assume without loss of generality that $\sum_{\gamma\in\Gamma}\alpha_\gamma=1$.
Thus we can write $x^{p}=\mathbb{E}[\boldsymbol{\gamma}_{1}\cdots\boldsymbol{\gamma}_{p}]$ for every $p\in\N$, where $\boldsymbol{\gamma}_{i}$ are i.i.d.\ random elements of $\Gamma$ with $\mathbb{P}\left[\boldsymbol{\gamma}_{i}=\gamma\right] = \alpha_{\gamma}$.
In particular, \eqref{eq:u1basic} yields
\[
u_{1}\left(x^{p}\right)= - \tau{\left(\lambda(x)^p\right)} +
\sum_{k=2}^{p\ell}\left(\omega(k)-1\right)\sum_{v\in\Gamma_{{\rm np}}}\mathbb{E}\left[1_{\boldsymbol{\gamma}_{1}\cdots\boldsymbol{\gamma}_{p}\,\equiv_{\Gamma}\,v^{k}}\right],
\]
where $\Gamma_{\rm np}$ denotes the set of non-powers in $\Gamma$.
Here we used that $\boldsymbol{\gamma}_{1}\cdots\boldsymbol{\gamma}_{p}\,\equiv_{\Gamma}\,v^{k}$
implies that $k\le|v^{k}|\le|\boldsymbol{\gamma}_{1}|+\cdots+|\boldsymbol{\gamma}_{p}|\le p\ell$ using Lemma~\ref{lem:linear_lower_bound} below,
and that $\E[1_{\boldsymbol{\gamma}_{1}\cdots\boldsymbol{\gamma}_{p}\,\equiv_{\Gamma}\,1}]=
\E[\tau(\lambda(\boldsymbol{\gamma}_{1}\cdots\boldsymbol{\gamma}_{p}))]=\tau(\lambda(x)^p)$.

By Lemma \ref{lem:padding power words to get nice structure}, we
have for every $k\ge2$ and $\gamma_{1},\ldots,\gamma_{p}\in\Gamma$
\begin{eqnarray*}
\sum_{v\in\Gamma_{\text{np}}}\mbox{\ensuremath{1_{\gamma_{1}\cdots\gamma_{p}\,\equiv_{\Gamma}\,v^{k}}}} & \le & \sum_{1\le t_{1}\le\cdots\le t_{4}\le p+1}\sum_{\substack{b,h,u_{1},\ldots,u_{4}\in\Gamma\\
|u_{i}|\le c_{3}(\ell+\log p)
}
}1_{\gamma_{1}\cdots\gamma_{t_{1}-1}\,\equiv_{\Gamma}\,bu_{1}^{-1}}1_{\gamma_{t_{1}}\cdots\gamma_{t_{2}-1}\,\equiv_{\Gamma}\,u_{1}hu_{2}^{-1}}\times\\
 &  & \,1_{\gamma_{t_{2}}\cdots\gamma_{t_{3}-1}\,\equiv_{\Gamma}\,u_{2}hu_{3}^{-1}}\,1_{\gamma_{t_{3}}\cdots\gamma_{t_{4}-1}\,\equiv_{\Gamma}\,u_{3}h^{k-2}u_{4}^{-1}}\,1_{\gamma_{t_{4}}\cdots\gamma_{p}\,\equiv_{\Gamma}\,u_{4}b^{-1}}.
\end{eqnarray*}
As 
\[
1_{\gamma_{1}\cdots\gamma_{p}\,\equiv_{\Gamma}\,\gamma}=\left\langle \delta_{\gamma},\lambda(\gamma_{1}\cdots\gamma_{p})\delta_{e}\right\rangle ,
\]
this yields 
\begin{multline*}
\sum_{v\in\Gamma_{\text{np}}}\mbox{\ensuremath{\mathbb{E}\left[1_{\boldsymbol{\gamma}_{1}\cdots\boldsymbol{\gamma}_{p}\,\equiv_{\Gamma}\,v^{k}}\right]}}  \le  \sum_{1\le t_{1}\le\cdots\le t_{4}\le p+1}\sum_{\substack{b,h,u_{1},\ldots,u_{4}\in\Gamma\\
|u_{i}|\le c_{3}(\ell+\log p)
}
}\left\langle \delta_{bu_{1}^{-1}},\lambda(x)^{t_{1}-1}\delta_{e}\right\rangle \times
\\
 \left\langle \delta_{u_{1}hu_{2}^{-1}},\lambda(x)^{t_{2}-t_{1}}\delta_{e}\right\rangle ~
\left\langle \delta_{u_{2}hu_{3}^{-1}},\lambda(x)^{t_{3}-t_{2}}\delta_{e}\right\rangle ~ 
\left\Vert \lambda(x)\right\Vert ^{t_{4}-t_{3}} ~
\left\langle \delta_{u_{4}b^{-1}},\lambda(x)^{p+1-t_{4}}\delta_{e}\right\rangle.
\end{multline*}
Now note that the sum of $t_{1},\ldots,t_{4}$ has at most $(p+1)^{4}$
terms, while the sum over $u_{1},\ldots,u_{4}$ has at most $(4g+1)^{4c_{3}(\ell+\log p)}$
terms. Moreover, as 
\[
\sum_{v\in\Gamma}|\langle\delta_{v},\lambda(x)^{t}\delta_{e}\rangle|^{2}=\|\lambda(x)^{t}\delta_{e}\|^{2}\le\|\lambda(x)\|^{2t},
\]
we can apply Cauchy-Schwarz to the remaining sums over $b,h$ to estimate
\begin{align*}
u_{1}(x^{p}) &\le |\tau{\left(\lambda(x)^p\right)}|+(p\ell)^{2}\max_{2\le k\le p\ell}\sum_{v\in\Gamma_{{\rm np}}}\mathbb{E}\big[1_{\boldsymbol{\gamma}_{1}\cdots\boldsymbol{\gamma}_{p}\,\equiv_{\Gamma}\, v^{k}}\big]
\\
&\le\left[1+(p\ell)^{2}(p+1)^{4}(4g+1)^{4c_{3}(\ell+\log p)}\right]\|\lambda(x)\|^{p}.
\end{align*}
The conclusion follows directly. 
\end{proof}

Above we used the following lemma, whose proof we spell out for completeness. An analogous property holds up to a constant for any torsion-free hyperbolic group by \cite[Proposition~III.$\Gamma$.3.15 and Corollary~III.$\Gamma$.3.10]{bridson2013metric}.

\begin{lem}
\label{lem:linear_lower_bound}
For every $1\ne v\in\Gamma$ and $k\in\Z_+$, we have $k\le|v^k|$.
\end{lem}

\begin{proof}
Write $v=b\gamma b^{-1}$ where $\gamma$ has minimal word length in the conjugacy class of $v$. Then $|\gamma^n|=n|\gamma|$ for all
$n\in\Z_+$ by \cite[Lemma~2.11]{BirmanSeries}. Therefore
$nk|\gamma| = |\gamma^{nk}| = |b^{-1}v^{nk}b| \le n|v^k| + 2|b|$ for all $k,n\in\Z_+$, and the
conclusion follows by letting $n\to\infty$.
\end{proof}

\section{Proof of Theorem \ref{thm:strong_convergence}\label{sec:Main-proof}}

Before we complete the proof, we must dispense
with a cosmetic issue. Let
\[
\pi_{n}=\mathrm{std}\circ\phi_{n},
\]
where $\phi_{n}$ is a uniformly distributed random element of $\mathrm{Hom}(\Gamma_{g},S_{n})$. In Theorem \ref{thm:asymp-effective}, we have proved an effective
asymptotic expansion of $\E[\tr \pi_n(\gamma)]$ in powers of $\frac{1}{n}$. However, since
$\pi_n:\Gamma\to\U(n-1)$ is a representation of dimension $n-1$, this differs slightly
from the setting of \S\ref{subsec:Technical-innovation-I:surpassing-poly} where $\pi_n$ is $n$-dimensional.

This minor discrepancy is readily resolved by noting that the statement and proof of Theorem~\ref{thm:main} extend \emph{verbatim} to the case that $\pi_n$ is $(n-1)$-dimensional, provided that we still understand the normalized trace to be defined by 
\[
    \ntr \pi_n \eqdf \frac{1}{n}\tr \pi_n.
\]
We will therefore impose this convention in the remainder of this section.\footnote{The issue is purely cosmetic in nature, since the proof of Theorem \ref{thm:asymp-effective} could also be readily adapted to yield an effective asymptotic expansion in powers of $\frac{1}{n-1}$.}

We can now proceed to assemble all the pieces of the proof. For completeness, we first 
spell out a quantitative strong convergence bound.

\begin{thm}
\label{thm:strongquant}
For any $x\in \mathrm{M}_d(\C)\otimes\C[\Gamma_g]$, $n\ge1$, and $\varepsilon>0$,
we have
\[
\P\big[\|[\mathrm{id}\otimes\pi_{n}](x)\|\ge(1+\varepsilon)\|[\mathrm{id}\otimes\lambda](x)\|\big]\le\frac{cd}{n\varepsilon^b}.
\]
Here $b$ is a constant that depends only on $g$, and $c$ depends on
$g$ and $|x|$.
\end{thm}

\begin{proof}
Assume first that $x$ is self-adjoint.
The validity of the effective asymptotic expansion of Assumption \ref{ass:cond1} follows immediately from Theorem \ref{thm:asymp-effective} and \eqref{eq:fixrep} with
\[
    u_1(\gamma) = \a_0(\gamma)-1,\qquad\quad
    u_k(\gamma) = \a_{k-1}(\gamma)\quad\text{for}\quad k\ne 1.
\]
Moreover, it is shown in \cite[Thm.~1.2]{MPasympcover} that $\a_{-1}(\gamma)=1_{\gamma=1}$ and
$\a_{0}(\gamma)=\omega(\gamma)$. Thus
\[
    u_0(\gamma) = 1_{\gamma=1} = \tau(\lambda(\gamma)),\qquad\qquad 
    u_1(\gamma) = \omega(\gamma)-1.
\]
Thus Assumption \ref{ass:cond1} holds, and Assumption \ref{ass:cond2bis} follows from
Theorem \ref{thm:completepowers}. The conclusion therefore follows from Theorem \ref{thm:quantmtx} in the case that $x$ is self-adjoint.
If $x$ is not self-adjoint, the conclusion follows by applying the self-adjoint bound to $x\leftarrow x^*x$.
\end{proof}

We finally complete the proof of Theorem \ref{thm:strong_convergence}.

\begin{proof}[Proof of Theorem \ref{thm:strong_convergence}]
Theorem \ref{thm:strongquant} immediately yields the upper bound
\[
    \|\pi_n(x)\| \le \|\lambda(x)\| + o(1) \quad\text{with probability}\quad 1-o(1)
\]
as $n\to\infty$ for every $x\in\C[\Gamma]$. 
As surface groups have the unique trace property \cite{delaharpe88},
the corresponding strong convergence lower bound
\[
    \|\pi_n(x)\| \ge \|\lambda(x)\| - o(1) \quad\text{with probability}\quad 1-o(1)
\]
follows from the upper bound, see, e.g., \cite[\S 5.3]{MageeSalle2}.
\end{proof}

\bibliographystyle{amsalpha}
\bibliography{sc_random_surfaces}

~\\
Michael Magee\\
Department of Mathematical Sciences, Durham University, DH1 3LE Durham, UK\\
\noindent\texttt{michael.r.magee@durham.ac.uk}\\
\\
Doron Puder\\
School of Mathematical Sciences, Tel Aviv University, Tel Aviv, 6997801, Israel\\
and IAS Princeton, School of Mathematics, 1 Einstein Drive, Princeton 08540, USA\\
\texttt{doronpuder@gmail.com}~\\
\\
Ramon van Handel\\
Department of Mathematics, Princeton University, Princeton, NJ 08544, USA\\
\texttt{rvan@math.princeton.edu}

\end{document}